\documentclass[a4paper,reqno, oneside]{amsart}
\usepackage[english]{babel}
\usepackage{amsmath, amssymb, amsthm, amscd}
\usepackage{enumerate}
\usepackage{palatino}
\usepackage{mathpazo}
\usepackage{paralist}
\usepackage[a4paper]{geometry}

\usepackage{hyperref}

 \usepackage{amsaddr}

\usepackage{tensor}

\usepackage{appendix}

\usepackage{tikz} 
\usepackage{tikz-cd}

\usepackage{bm}

\usepackage{mathtools}

\DeclareMathOperator{\id}{id}

\makeatletter
\newcommand{\ostar}{\mathbin{\mathpalette\make@circled\ast}}
\newcommand{\make@circled}[2]{%
  \ooalign{$\m@th#1\smallbigcirc{#1}$\cr\hidewidth$\m@th#1#2$\hidewidth\cr}%
}
\newcommand{\smallbigcirc}[1]{%
  \vcenter{\hbox{\scalebox{0.77778}{$\m@th#1\bigcirc$}}}%
}
\makeatother

\renewcommand{\setminus}{\smallsetminus}

\theoremstyle{plain}
\newtheorem{theorem}{Theorem}[section]
\newtheorem{prop}[theorem]{Proposition}
\newtheorem{lemma}[theorem]{Lemma}
\newtheorem{cor}[theorem]{Corollary}  

\newtheorem*{thm}{Theorem}

\theoremstyle{definition}
\newtheorem{definition}[theorem]{Definition}

\theoremstyle{remark}
\newtheorem{remark}[theorem]{Remark}
\newtheorem{example}[theorem]{Example}

\numberwithin{equation}{section}

\begin{document}
\setlength{\parindent}{0.cm}

\title{String topology on the space of paths with endpoints in a submanifold}

\author{Maximilian Stegemeyer}
\address{Mathematisches Institut, Universit\"at Freiburg, Ernst-Zermelo-Straße 1, 79104 Freiburg, Germany}
\email{maximilian.stegemeyer@math.uni-freiburg.de}
\date{\today}
\keywords{String Topology, Loop Spaces}
\subjclass{55P50}

\begin{abstract}
In this article we consider algebraic structures on the homology of the space of paths in a manifold with endpoints in a submanifold.
The Pontryagin-Chas-Sullivan product on the homology of this space had already been investigated by Hingston and Oancea for a particular example. 
We consider this product as a special case of a more general construction where we consider pullbacks of the path space of a manifold under arbitrary maps.
The product on the homology of this space as well as the module structure over the Chas-Sullivan ring are shown to be invariant under homotopies of the respective maps.
This in particular implies that the Pontryagin-Chas-Sullivan product as well as the module structure on the space of paths with endpoints in a submanifold are isomorphic for two homotopic embeddings of the submanifold.
Moreover, for null-homotopic embeddings of the submanifold this yields nice formulas which we can be used to compute the product and the module structure explicitly.
We show that in the case of a null-homotopic embedding the homology of the space of paths with endpoints in a submanifold is even an algebra over the Chas-Sullivan ring.
\end{abstract}

\maketitle

\tableofcontents

\section{Introduction}

String topology is usually understood as the study of algebraic structures on the homology or cohomology of the free loop space of a closed oriented manifold.
The most prominent example of a string topology operation is the Chas-Sullivan product which takes the form
$$    \wedge_{\mathrm{CS}} \colon \mathrm{H}_i(\Lambda M)\otimes \mathrm{H}_j(\Lambda M)\to \mathrm{H}_{i+j-n}(\Lambda M).    $$
Here, $\Lambda M$ is the free loop space of the closed oriented manifold $M$ and $n$ is the dimension of $M$.
The Chas-Sullivan product is associative, unital and graded commutative.
The geometric idea behind the Chas-Sullivan product is to concatenate those loops in the image of two cycles which share the same basepoint.
It is therefore related to the intersection product on the base manifold $M$.
We refer to \cite{chas:1999} for the first definition of the Chas-Sullivan product and to \cite{cohen:2002} and \cite{hingston:2017} for equivalent definitions of the Chas-Sullivan product which are closer to the one which is used in this manuscript.

In the present article we consider the space of paths in a manifold with endpoints in a given submanifold.
More precisely, let $M$ be a closed oriented manifold and let $N$ be a closed oriented submanifold.
Set
$$      P_N M = \{ \gamma\in PM\,|\, \gamma(0),\gamma(1)\in N\}     $$
where $PM$ is the space of paths $\gamma\colon [0,1]\to M$.
If we need to specify the embedding $i\colon N\to M$ we shall also write $P_N^i M$.
The homology of the path space $P_N M$ has been studied in \cite{abbondandolo:2008} where the authors prove that the Floer homology for Hamiltonian orbits on the cotangent bundle of a compact
manifold with certain non-local conormal boundary conditions is isomorphic to the singular homology of a path space of the type $P_N M$.
String topological ideas on the path space $P_N M$ have been examined in \cite{blumberg:2009} where the authors study open-closed topological conformal field theories.

In this paper, we will study the \textit{path product}
$$      \wedge\colon \mathrm{H}_i(P_N M)\otimes \mathrm{H}_j(P_N M) \to \mathrm{H}_{i+j-k}(P_N M)     $$
where $k$ is the dimension of the submanifold $N$.
Motivated by the Chas-Sullivan product the idea of its definition is to use the concatenation of paths $\gamma,\sigma\in P_N M$ satisfying $\gamma(1) = \sigma(0)$.
It is clear that this definition is not symmetric, so we cannot expect to get a commutative product.
The path product is however associative and unital.

As far as the author is aware this product has only been studied so far by Hingston and Oancea in \cite{hingston:2013oancea} where the authors study the particular case $M = \mathbb{C}P^n$ and $N = \mathbb{R}P^n$.
Hingston and Oancea use a Morse-theoretic description of the homology of the path space $P_{\mathbb{R}P^n} \mathbb{C}P^n$ to compute the product.

In addition to the path product we shall also consider the following module structure.
There is a pairing
\begin{equation}\label{eq_module_introduction}
        \cdot \colon \mathrm{H}_i(\Lambda M)\otimes \mathrm{H}_j(P_N M)\to \mathrm{H}_{i+j-n}(P_N M)    
\end{equation}
which makes $\mathrm{H}_{\bullet}(P_N M)$ into a left module over the Chas-Sullivan ring $(\mathrm{H}_{\bullet}(\Lambda M),\wedge_{\mathrm{CS}})$.
It is then an obvious question if the homology of $P_N M$ together with the path product $\wedge$ is an algebra over the Chas-Sullivan ring.
The heuristic geometric expectation is that this should not be the case, since the definition of the path product is too asymmetric.
Moreover, it is an important question to understand whether the path product and the module structure are invariant under homotopies of the embedding $i\colon N\to M$ of the submanifold.

In the present article we study these questions by considering a more general construction of the path product and the module structure which we shall briefly introduce now.
Let $M$ and $N$ be closed oriented manifolds and let $f\colon N\to M$ be a smooth map, not necessarily an embedding.
Then we take the pullback of the free path fibration $PM\to M\times M,\quad \gamma\mapsto (\gamma(0),\gamma(1))$
under the map $(f\times f)\colon N\times N\to M\times M$.
The pullback will be denoted by $P^f$ and it is explicitly given as
$$   P^f =  \{ (x_0,x_1,\gamma)\in N\times N\times PM\,|\,  \gamma(0) = f(x_0),\,\, \gamma(1) = f(x_1) \} .     $$
If $f\colon N\to M$ is the inclusion of a submanifold then we recover the space $P_N M$, i.e. the space of paths in $M$ with endpoints in $N$.
In this article we define a product 
$$   \wedge\colon \mathrm{H}_i(P^f) \otimes \mathrm{H}_j(P^f) \to \mathrm{H}_{i+j-k}(P^f)       $$
where $k = \mathrm{dim}(N)$ which agrees with the path product that we explained above in the case that $f\colon N\to M$ is the inclusion of a submanifold.
The advantage of considering this more general construction is that it is easier to understand how the product behaves for maps homotopic to $f$, where we allow arbitrary maps and not just embeddings of submanifolds.
We shall also generalize the module structure on $\mathrm{H}_{\bullet}(P_N M)$ over the Chas-Sullivan ring by introducing a module structure
$$    \cdot \colon \mathrm{H}_i(\Lambda M)\otimes \mathrm{H}_j(P^f) \to \mathrm{H}_{i+j-n}(P^f) .  $$
Again, this agrees with the pairing in equation \eqref{eq_module_introduction} in case that $f$ is the embedding of a submanifold.

Let us now summarize the main results of this paper.
We shall study how the product and the module structure on $\mathrm{H}_{\bullet}(P^f)$ behave if we compare a map $f\colon N\to M$ with a map $g\colon N\to M$ which is homotopic to $f$.
The spaces $P^f$ and $P^g$ are then homotopy equivalent.
We show that the induced isomorphism in homology also induces isomorphisms of the ring structure and of the module structure.
\begin{thm}[Theorem \ref{theorem_invariance_path_product} and Theorem \ref{theorem_invariance_module_structure}]
    Let $M$ and $N$ be closed oriented manifolds and assume that $f\colon N\to M$ and $g\colon N\to M$ are smooth maps.
    Take homology with coefficients in a field $R$ and consider the homology of the spaces $P^f$ and $P^g$.
    If $f$ and $g$ are homotopic then
    \begin{enumerate}
        \item the rings $(\mathrm{H}_{\bullet}(P^f), \wedge)$ and $(\mathrm{H}_{\bullet}(P^g), \wedge)$ are isomorphic and
        \item the homologies $\mathrm{H}_{\bullet}(P^f)$ and $\mathrm{H}_{\bullet}(P^g)$ are isomorphic as left modules over the Chas-Sullivan ring of $M$.
    \end{enumerate}
\end{thm}

It is a direct consequence of this theorem that if $f\colon N\hookrightarrow M$ and $g\colon N\hookrightarrow M$ are two homotopic embeddings of $N$ as a submanifold of $M$, then the ring and the module structure on $\mathrm{H}_{\bullet}(P_N^f M)$ and $\mathrm{H}_{\bullet}(P_N^g M)$ are isomorphic.

Moreover, the above theorem can be used to study the following situation.
If the embedding $f\colon N\hookrightarrow M$ is null-homotopic, then one sees that
$$    P_N M \simeq N\times N\times \Omega M     $$
where $\Omega M$ is the based loop space of $ M$.
It is now an obvious question if we can describe the path product and the module structure intrinsically on the tensor product $\mathrm{H}_{\bullet}(N)^{\otimes 2}\otimes \mathrm{H}_{\bullet}(\Omega M)$.
Take homology with coefficients in a field $\mathbb{K}$ and define
$$    \mu\colon  \mathrm{H}_{\bullet}(N)^{\otimes 2}\otimes \mathrm{H}_{\bullet}(N)^{\otimes 2} \to  \mathrm{H}_{\bullet}(N)^{\otimes 2}  $$
by
$$     \mu(a\otimes b\otimes c\otimes d) =   h_* (b\overline{\cap}c) a\otimes d   \quad \text{for}\,\,\,a,b,c,d\in\mathrm{H}_{\bullet}(N)  $$
where $\overline{\cap}$ is the intersection product of $N$ and where $h\colon N\to \{\mathrm{pt}\}$ is the constant map.
We thus view $h_*$ as the map $h_*\colon \mathrm{H}_{\bullet}(N)\to \mathbb{K}$.
We equip the tensor product $\mathrm{H}_{\bullet}(N)^{\otimes 2}\otimes \mathrm{H}_{\bullet}(\Omega M)$ with the product
$$    \mu_{N,\Omega}\colon  \big(\mathrm{H}_{\bullet}(N)^{\otimes 2}\otimes \mathrm{H}_{\bullet}(\Omega M)\big) \otimes \big(\mathrm{H}_{\bullet}(N)^{\otimes 2} \otimes \mathrm{H}_{\bullet}(\Omega M)\big) \to  \mathrm{H}_{\bullet}(N)^{\otimes 2} \otimes \mathrm{H}_{\bullet}(\Omega M) $$
defined by
$$     \mu_{N,\Omega}((a\otimes b\otimes x)\otimes ( c\otimes d\otimes y)) =   (-1)^{|x|(|c|+|d|+k)} \mu(a\otimes b\otimes c\otimes d)\otimes (x\star y) 
$$
for
$ a,b,c,d\in\mathrm{H}_{\bullet}(N),\,\,x,y\in\mathrm{H}_{\bullet}(\Omega M)  $
and with $\star\colon \mathrm{H}_{\bullet}(\Omega M)\otimes \mathrm{H}_{\bullet}(\Omega M)\to \mathrm{H}_{\bullet}(\Omega M)$ being the Pontryagin product.
Furthermore, define a pairing
$$ \nu_{N,\Omega}\colon \mathrm{H}_{\bullet}(\Lambda M ) \otimes \big( \mathrm{H}_{\bullet}(N)^{\otimes 2}\otimes \mathrm{H}_{\bullet}(\Omega M)\big) \to \mathrm{H}_{\bullet}(N)^{\otimes 2}\otimes \mathrm{H}_{\bullet}(\Omega M)    $$
by setting
$$    \nu_{N,\Omega}( A\otimes a\otimes b\otimes x) =  (-1)^{(n+|A|)(n+|a|+|b|)} a\otimes b\otimes (j_! A \star x)   $$
for $A\in\mathrm{H}_{\bullet}(\Lambda M ),  a,b\in\mathrm{H}_{\bullet}(N),x\in\mathrm{H}_{\bullet}(\Omega M)$.
Here, $j_!\colon \mathrm{H}_i(\Lambda M)\to \mathrm{H}_{i-n}(\Omega M)$ is the Gysin map induced by the inclusion $j\colon \Omega M\to \Lambda M$.
We then show the following.
\begin{thm}[Proposition \ref{prop_trivial_map_path_space} and Proposition \ref{prop_trivial_map_module}]
    Let $M$ be a closed oriented manifold and $N$ a closed oriented submanifold with null-homotopic inclusion.
    Take homology with coefficients in a field.
    Then
    \begin{enumerate}
        \item the ring $(\mathrm{H}_{\bullet}(P_N M),\wedge)$ is isomorphic to the ring 
        $   (\mathrm{H}_{\bullet}(N)^{\otimes 2}\otimes \mathrm{H}_{\bullet}(\Omega M),\mu_{N,\Omega})     $ and
        \item the homology $\mathrm{H}_{\bullet}(P_N M)$ is isomorphic as a module over the Chas-Sullivan ring of $M$ to the tensor product $\mathrm{H}_{\bullet}(N)^{\otimes 2}\otimes \mathrm{H}_{\bullet}(\Omega M)$ equipped with the module structure $\mu_{N,\Omega}$.
    \end{enumerate}
\end{thm}

Note that if $M = \mathbb{S}^n$ is a sphere then any submanifold $N\hookrightarrow \mathbb{S}^n$ has null-homotopic inclusion.
Hence, the above theorem gives a complete description of the path product and the module structure in the case that $M$ is a sphere.
With the explicit descriptions of the previous theorem we show the following.
\begin{thm}[Theorem \ref{theorem_algebra_over_cs}]
    Let $M$ be a closed oriented manifold and $N$ a closed oriented submanifold with null-homotopic inclusion.
    Assume that the image of the Gysin map $j_!\colon \mathrm{H}_i(\Lambda M)\to \mathrm{H}_{i-n}(\Omega M)$ lies in the center of the Pontryagin ring $(\mathrm{H}_{\bullet}(\Omega M),\star)$.
    Then the homology $\mathrm{H}_{\bullet}(P_N M)$ equipped with the path product is an algebra over the Chas-Sullivan ring of $M$.
\end{thm}
We show that the condition on the Gysin map is satisfied if $M = \mathbb{S}^n$ is a sphere.
Note that this condition always holds if the Pontryagin ring is graded commutative.

This article is organized as follows.
In Section \ref{sec_def_path_product} we introduce the path product via the more general construction of the space $P^f$. We will then discuss the basic properties of the path product.
Certain classes of examples of the path product are investigated in Section \ref{sec_path_product_examples}.
In Section \ref{sec_invariance_product} we show that the product on $\mathrm{H}_{\bullet}(P^f)$ is invariant under homotopies of $f$.
The module structure over the Chas-Sullivan ring is introduced in Section \ref{sec_definition_module} and we study the behavior of the module structure in some specific cases in Section \ref{sec_module_structure_explicit}.
Finally, in Section \ref{sec_invariance_module} we show that the module structure on $\mathrm{H}_{\bullet}(P^f)$ is invariant under homotopies of $f$ and discuss consequences of this invariance.

\medskip

\emph{In this article we assume that all manifolds are connected. Moreover, if $X$ and $Y$ are topological spaces, the tensor product of homology groups $\mathrm{H}_{\bullet}(X) \otimes \mathrm{H}_{\bullet}(Y)$ is always understood as the tensor product over the coefficient ring, i.e. $\mathrm{H}_{\bullet}(X) \otimes \mathrm{H}_{\bullet}(Y) = \mathrm{H}_{\bullet}(X;R) \otimes_R \mathrm{H}_{\bullet}(Y;R)$.
If the coefficient ring $R$ is a field then we shall make the identification given by the Künneth isomorphism $\mathrm{H}_{\bullet}(X\times Y) \cong \mathrm{H}_{\bullet}(X)\otimes \mathrm{H}_{\bullet}(Y)$ without always explicitly stating it.}

\medskip

\noindent \textbf{Acknowledgements:} 
The author thanks the anonymous referee for many helpful comments and suggestions.

The revised version of this manuscript owes its existence to a discussion that the author had with Lie Fu and Alexandru Oancea during a visit at the IRMA Strasbourg in January 2024.
The author is very grateful for the hospitality of the IRMA.

Some of the work on this manuscript was carried out during a stay of the author at the Copenhagen Centre for Geometry and Topology supported by the Deutsche Forschungsgemeinschaft (German Research Foundation) – grant agreement number 518920559.
The author is grateful for the support by the Danish National Research Foundation through the Copenhagen Centre for Geometry and Topology (DNRF151).

\section{Definition and basic properties of the path product}\label{sec_def_path_product}

In this section we introduce the path space $P^f$ with respect to a map $f\colon N\to M$ between manifolds $N$ and $M$.
We then define the path product on the homology of $P^f$ and study its basic algebraic properties.

Let $M$ be a closed manifold and $N$ a closed oriented manifold.
Let $f\colon N\to M$ be a smooth map.
Let $R$ be a commutative ring and take homology with coefficients in $R$.
Moreover, let $n = \mathrm{dim}(M)$ and $k = \mathrm{dim}(N)$.
We consider the path space 
$$  PM = \Big\{ \gamma\colon [0,1]\to M\,|\, \gamma \text{ absolutely continuous }, \,\, \int_0^1|\Dot{\gamma}(s)|^2\,\mathrm{d}s <\infty \Big\}   ,  $$
with the notion of absolutely continuous curves in a manifold as defined in \cite{klingenberg:1995}.
In the following we shall denote the unit interval $[0,1]$ by $I$.
It turns out that the space $PM$ is a Hilbert manifold, see \cite[Section 2.3]{klingenberg:1995}.
The evaluation map 
\begin{equation} \label{eq_ev}
    \mathrm{ev}\colon PM\to M\times M,\quad \gamma\mapsto (\gamma(0),\gamma(1))
\end{equation} 
is both a submersion and a fibration.
Consider the pullback
$$
\begin{tikzcd}
    P^f \arrow[]{r}{} \arrow[swap]{d}{(p_0,p_1)} & [2.5em] PM \arrow[]{d}{} \\ N\times N \arrow[]{r}{f\times f} & M\times M .
\end{tikzcd}
$$
Explicitly, we have that
$$    P^f =  \{ ( x_0,x_1,\gamma)\in N\times N\times PM\,|\, \gamma(0) = f(x_0), \,\gamma(1) = f(x_1) \}   .   $$
By construction, the map
$$ p= (p_0,p_1)\colon P^f \to N\times N ,\quad (x_0,x_1,\gamma) \mapsto (x_0,x_1) $$
is a fibration.
Moreover since $\mathrm{ev}\colon PM\to M\times M$ is a submersion and since $f$ is smooth, we see that $P^f$ is a Hilbert manifold and that the map $p$ is a submersion.

\begin{remark}
    Our main interest on the spaces of the form $P^f$ is the case when $f\colon N\to M$ is the embedding of a submanifold.
    In this case we have a canonical homeomorphism
    $$  P^f \xrightarrow[]{\cong} P_N M : = \{ \gamma\in PM\,|\, \gamma(0),\gamma(1) \in N\}       $$
    given by mapping $(x_0,x_1,\gamma)\in P^f$ to $\gamma\in P_N M$.
    Therefore we shall also write $P^f = P_N M$ without explicitly stating this homeomorphism.
    If we want to specify the embedding $f$ we also write $P_N^f M$.
\end{remark}

We now turn to describing the path product $\wedge\colon \mathrm{H}_i(P^f)\otimes \mathrm{H}_j(P^f) \to \mathrm{H}_{i+j-k}(P^f)$.
Let $\Delta N\subseteq N\times N$ be the diagonal.
Since the maps $p_0,p_1\colon P^f \to N$ given by
$$    p_i(x_0,x_1,\gamma) =  x_i\quad \text{for}\quad i\in \{0,1\}      $$
are fibrations and a submersions we have that
$$    C^f =   (p_1\times p_0)^{-1}(\Delta N) = \{ ((x_0,x_1,\gamma),(y_0,y_1,\sigma))\in P^f\times P^f\,|\, x_1 = y_0\}  $$
is a codimension $k$-submanifold of $P^f\times P^f$ where $k = \mathrm{dim}(N)$.
In particular for $(\gamma,\sigma)\in C^f$ we have $\gamma(1) = \sigma(0)$.
There is a map $C^f\to P^f$ given by
$$     c\colon ((x_0,z,\gamma),(z,y_1,\sigma)) =  (x_0,y_1,  \mathrm{concat}(  \gamma , \sigma)) .   $$
Here, $\mathrm{concat}$ is the concatenation of two paths for which we make the following convention.
Define $F_2\subseteq PM \times PM$ to be the space
$$   F_2 = \{(\gamma,\sigma)\in PM\times PM\,|\, \gamma(1) = \sigma(0) \} .     $$
Then we let $\mathrm{concat}\colon F_2\to PM$ be the map
$$\mathrm{concat}(\gamma,\sigma)(t) = \begin{cases}
                \gamma(2t), & 0\leq t\leq\tfrac{1}{2} \\ \sigma(2t-1), & \tfrac{1}{2}\leq t \leq 1.
\end{cases}     $$
For $m\geq 2$ there are corresponding spaces $F_m$ and maps $\mathrm{concat}\colon F_m\to PM$ given by
$$  \mathrm{concat}(\gamma_1,\ldots,\gamma_m) = \begin{cases}
                \gamma_1(mt), & 0\leq t\leq\tfrac{1}{m} \\
                \vdots, & \vdots
                \\ \gamma_m(mt-(m-1)) , & \tfrac{m-1}{m}\leq t\leq 1 .
\end{cases}     $$  
In this article we will frequently use these concatenation maps and their restriction to certain subspaces of $F_m$.
It will be clear from the number of arguments and the context how the respective concatenation map is understood.

We now describe a tubular neighborhood of the diagonal $\Delta N\hookrightarrow N\times N$.
Fix a Riemannian metric $g$ on $M$ and let $\epsilon'>0$ be smaller than the injectivity radius on $M$.
Furthermore choose a Riemannian metric $h$ on $N$ and an $\epsilon> 0$ smaller than the injectivity radius of $N$.
By compactness of $N$ we can choose $\epsilon$ such that $\mathrm{d}_N(p,q)<\epsilon$ implies $\mathrm{d}_M(f(p),f(q))<\epsilon'$.
Here, $\mathrm{d}_N$ is the distance function on $N$ induced by $h$ and $\mathrm{d}_M$ is the distance function on $M$ induced by $g$.
Then we set 
$$   U_N = \{ (p,q)\in N\times N\,|\,  \mathrm{d}_N(p,q)<\epsilon   \}     $$
and one can check that this is indeed a tubular neighborhood of the diagonal, see \cite{hingston:2017}.
In particular we have that
$$   \mathrm{H}^k(U_N,U_N\setminus \Delta N)\cong \mathrm{H}^k(TN, TN\setminus N)    $$
hence there is a Thom class $\tau_N\in \mathrm{H}^k(U_N,U_N\setminus \Delta N)$.

Now, consider the map of pairs $p_1\times p_0\colon P^f\times P^f\to N\times N$ and the preimages
$$    U_{C^f} =    (p_1\times p_0)^{-1} (U_N)   =   \{((x_0,x_1,\gamma),(y_0,y_1,\sigma)\in P^f\times P^f\,|\, (x_1,y_0)\in U_N\} .      $$
We pull back the Thom class to get a class $\tau_{C^f} = (p_1\times p_0)^* \tau_N \in \mathrm{H}^k(U_{C^f},U_{C^f}\setminus C^f)$.
We now describe a homotopy retraction $\mathrm{R}_C\colon U_{C^f}\to C^f$.
If $((x_0,x_1,\gamma),(y_0,y_1,\sigma))\in U_{C^f}$ then by definition we have $\mathrm{d}_N(x_1,y_0)<\epsilon$.
Hence, there is a unique length-minimizing geodesic in $N$ connecting $x_1$ and $y_0$.  
We denote this geodesic by $\widehat{x_1 y_0}\colon I\to N$.
Define $\mathrm{R}_{C^f}\colon U_{C^f}\to C^f$ by setting
$$     \mathrm{R}_{C^f}((x_0,x_1,\gamma),(y_0,y_1,\sigma)) = ((x_0,x_1,\gamma),(x_1,y_1,\mathrm{concat}(f(\widehat{x_1 y_0}),\sigma))) .       $$
We note here that one can indeed show explicitly that $U_{C^f}$ is a tubular neighborhood of $C_f$ and that the map $\mathrm{R}_{C^f}$ is homotopic to the retraction induced by the tubular neighborhood.
This can be done by adapting the strategy used in \cite[Proposition 2.2 and Lemma 2.3]{hingston:2017} to our current situation. 
In Section \ref{sec_definition_module} we encounter a very similar situation in which we prove the existence of a tubular neighborhood explicitly using these ideas.
Therefore we omit the proof that $U_{C^f}$ is a tubular neighborhood of $C_f$ here and refer to Lemma \ref{lemma_retraction_explicit} and \cite{hingston:2017}.
We can now define the product.

\begin{definition}\label{def_product}
    Let $M$ be a closed manifold and let $N$ be a closed oriented manifold of dimension $k$.
    Let $f\colon N\to M$ be a smooth map.
    Take homology with coefficients in a commutative ring $R$.
    Then the \textit{path product} of $P^f$ is defined as the composition
\begin{eqnarray*}
   \wedge \colon          
   \mathrm{H}_i(P^f)\otimes \mathrm{H}_j(P^f) 
   & \xrightarrow[]{\hphantom{}(-1)^{k-ki}\times\hphantom{}}  
   &  \mathrm{H}_{i+ j}(P^f\times P^f) \\
        & \xrightarrow[]{\hphantom{blaibalabla}} & \mathrm{H}_{i+j}(P^f\times P^f,P^f\times P^f \setminus C^f) 
        \\
        & \xrightarrow[]{\hphantom{il}\text{excision}\hphantom{l}} & \mathrm{H}_{i+j}(U_{C^f},U_{C^f}\setminus C^f)
        \\
            & \xrightarrow[]{\hphantom{blii}\tau_{C^f}\cap\hphantom{bl}} & \mathrm{H}_{i+j-k}(U_{C^f})
        \\
        & \xrightarrow[]{\hphantom{bi}(\mathrm{R}_{C^f})_*\hphantom{b}} & \mathrm{H}_{i+j-k}(C^f)
            \\ & \xrightarrow[]{\hphantom{blai}c_*\hphantom{blaii}} & \mathrm{H}_{i+j-k}(P^f) .    
\end{eqnarray*}
\end{definition}

\begin{remark}
    In the case that $f\colon N\to M$ is the embedding of a submanifold then we have seen that the space $P^f$ is the space of paths in $M$ with endpoints in the submanifold $N$.
    In this case this product can be seen to be the \textit{Pontryagin-Chas-Sullivan} product as described by Hingston and Oancea \cite{hingston:2013oancea}.
    Since we mostly deal with the more general situation of the path spaces $P^f$ and since we want to keep the exposition of this paper brief we decided to use the name \textit{path product} in this manuscript.
\end{remark}

In the following theorem we note some basic properties of the path product.
Before we do so we must introduce a product on the tensor product $\mathrm{H}_{\bullet}(N)\otimes \mathrm{H}_{\bullet}(N)$.
We begin with a general consideration.
Let $A$ be an $R$-module with $R$ a commutative ring.
Let $\beta\colon A\otimes A\to R$ be a bilinear form.
Then setting
$$   \mu_{\beta}\colon (A\otimes A)\otimes (A\otimes A) \to A\otimes A\quad \text{with}\quad \mu_{\beta}(a\otimes b\otimes c\otimes d) = \beta(b\otimes c)  a\otimes d $$
defines an associative product.
In our situation we will have $A = \mathrm{H}_{\bullet}(N)$.
Moreover, denote the intersection product on $N$ by $\overline{\cap}$ and let $1_N\in \mathrm{H}^0(N)$ be the unit of the cohomology ring of $N$.
Consider the bilinear form
$$  \beta\colon \mathrm{H}_{\bullet}(N)\otimes \mathrm{H}_{\bullet}(N)\to R, \quad \beta(X\otimes Y) = \begin{cases}
    \langle 1_N , X\overline{\cap} Y\rangle, & |X| + |Y| = k, \\ 0 & \text{else} .
\end{cases}    $$
Note that $\beta$ equals the bilinear form 
$$    \mathrm{H}_{\bullet}(N)\otimes \mathrm{H}_{\bullet}(N)\to R,\quad X\times Y\mapsto h_* (X\overline{\cap}Y)    $$
where $h\colon N\to \{\mathrm{pt}\}$ is the map to a point and thus $h_*\colon \mathrm{H}_{\bullet}(N)\to \mathrm{H}_0(\{\mathrm{pt}\})\cong R$.
We call this form the \textit{intersection product form} on $N$.
Thus by the above we obtain a product $$\mu_{\beta}\colon  (\mathrm{H}_{\bullet}(N)\otimes \mathrm{H}_{\bullet}(N))^{\otimes 2} \to  \mathrm{H}_{\bullet}(N)\otimes \mathrm{H}_{\bullet}(N)  $$
which is explicitly given by
$$    \mu_{\beta} (a\otimes b\otimes c\otimes d) =  h_* (b\overline{\cap} c) a\otimes d \quad\text{for} \,\,\, a,b,c,d\in\mathrm{H}_{\bullet}(N) .    $$
Moreover, we want to recall from \cite{hingston:2017} how the intersection product is defined in terms of a Gysin map.
Let $\tau_N\in\mathrm{H}^k(U_N,U_N\setminus \Delta N)$ be the Thom class as in the discussion before Definition \ref{def_product}.
Then the intersection product is the composition
\begin{eqnarray*}
   \overline{\cap} \colon            \mathrm{H}_i(N)\otimes \mathrm{H}_j(N ) & \xrightarrow[]{\hphantom{}(-1)^{k-ki}\times\hphantom{}}  &  \mathrm{H}_{i+ j}(N \times N ) \\
            &  \xrightarrow[]{\hphantom{blaiis_!}\hphantom{blika} } 
            & \mathrm{H}_{i+j}(U_N,U_N \setminus N) \\
            &  \xrightarrow[]{\hphantom{bail}\tau_N\cap\hphantom{bil}} & \mathrm{H}_{i+j-k}(U_N) 
              \\ 
              & \xrightarrow[]{\hphantom{bl}(\mathrm{R}_N)_*\hphantom{bil}} & \mathrm{H}_{i+j-k}(N) ,
\end{eqnarray*}
see \cite[Appendix B]{hingston:2017}.
Here, $\mathrm{R}_N\colon U_N\to N$ is the retraction of the tubular neighborhood.
Note that we use the non-standard notation $\overline{\cap}$ to distinguish the intersection product from the cap-product.
Without proof we note the following properties of the product $\mu_{\beta}$ which is induced by the intersection product form.
\begin{prop}\label{prop_intersection_pr_from}
    Let $N$ be a closed oriented manifold and take homology with coefficients in a field $\mathbb{K}$.
    Let $\beta\colon \mathrm{H}_{\bullet}(N)\otimes \mathrm{H}_{\bullet}(N)\to \mathbb{K}$ be the intersection product form and consider the induced product $\mu_{\beta}$.
    Then $\mu_{\beta}$ is an associative and unital product.
    If $N$ is of positive dimension then $\mu_{\beta}$ is not commutative.
\end{prop}

We now study the algebraic properties of the path product.

\begin{theorem}\label{theorem_properties_path_product}
    Let $M$ be a closed manifold and $N$ a closed oriented manifold with a smooth map $f\colon N\to M$.
     Take homology with coefficients in a commutative ring $R$.
    \begin{enumerate}
        \item The path product is an associative and unital product.
        \item Consider the fibration $p\colon P^f\to N\times N$ and assume that $R$ is a field.
        The induced map in homology $p_*\colon \mathrm{H}_{\bullet}(P^f)\to \mathrm{H}_{\bullet}(N)\otimes \mathrm{H}_{\bullet}(N)$ becomes a morphism of algebras if we equip $\mathrm{H}_{\bullet}(N)\otimes \mathrm{H}_{\bullet}(N)$ with the product $\mu_{\beta}$ induced by the intersection product form.
    \end{enumerate}
\end{theorem}
\begin{proof}
	The proof of the first part is very similar to the proof of the respective properties of the Chas-Sullivan product.
	Therefore we stay close to \cite[Appendix B]{hingston:2017}.
	Hingston and Wahl introduce a special relative version of the cap product, see \cite[Appendix A]{hingston:2017}, which we shall use in this proof.
	We need to introduce a thickened version the boundary of $U_{N}$ which we denote by
	$$   U_{N,\geq\epsilon_0} = \{ (p,q)\in U_N \,|\, \mathrm{d}(p,q)\geq \epsilon_0\}  .    $$
	Consider its preimage
	$$    U_{C^f,\geq\epsilon_0} = (p_1\times p_0)^{-1} (U_{N,\geq\epsilon_0}) .    $$
	The inclusion $(U_{C^f},U_{C^f,\geq\epsilon_0}) \hookrightarrow (U_{C^f}, U_{C^f}\setminus C^f)$ indcues an isomorphism in homology and therefore we can see the Thom class $\tau_{C^f}$ as a class $\tau_{C^f}\in\mathrm{H}^k(U_{C^f},U_{C^f,\geq\epsilon_0})$.
	The relative cap product is then defined to be a map
	$$   \mathrm{H}^k(U_{C^f},U_{C^f,\geq\epsilon_0})\otimes \mathrm{H}_i(P_{N}M) \to \mathrm{H}_{i-k}(U_{C^f}) .     $$
	See \cite[Appendix A]{hingston:2017} for details on this relative cap product.
	In the following commutative diagram we shall always use this cap product.
	For the rest of this proof we shall use the shorthand notation $P$ for $P^f$, $U$ for $U_{C^f}$ as well as $\mathrm{R}$ for the retraction $\mathrm{R}_{C^f}$ and $\tau_C$ for the Thom class $\tau_{C^f}$.
	$$
	\begin{tikzcd}[row sep=large]
		\mathrm{H}_i(P)\otimes \mathrm{H}_j(P)\otimes \mathrm{H}_l(P) \arrow[]{r}{\times} \arrow[]{d}{\times} &
		\mathrm{H}_{i+j}(P^2)\otimes \mathrm{H}_l(P) \arrow[]{r}{(\tau_C\cap)\otimes 1}\arrow[]{d}{\times} \arrow[dr, phantom,shift left=0.8ex, "(1)"] &
		\mathrm{H}_{i+j-k}(U)\otimes \mathrm{H}_l(P) \arrow[]{d}{\times}
		\\
		\mathrm{H}_i(P)\otimes \mathrm{H}_{j+l}(P^2) \arrow[]{r}{\times}\arrow[]{d}{1\otimes (\tau_C\cap)}\arrow[dr, phantom,shift left=0.8ex, "(2)"] &
		\mathrm{H}_{i+j+l}(P^3) \arrow[]{r}{(\tau_C\times 1)\cap}\arrow[]{d}{(1\times \tau_C)\cap} &
		\mathrm{H}_{m}(U\times P) \arrow[]{d}{(c \mathrm{R})_*\times 1}
		\\
		\mathrm{H}_i(P)\otimes \mathrm{H}_{j+l-k}(U) \arrow[]{r}{\times}\arrow[]{d}{1\otimes (c \mathrm{R})_*} &
		\mathrm{H}_{m}(P\times U) \arrow[]{d}{1\times (c \mathrm{R})_*}   \arrow[dr, phantom,shift left=0.8ex, "(3)"]  &
		\mathrm{H}_{m}(P\times P) \arrow[]{d}{\tau_C \cap} 
		\\
		\mathrm{H}_i(P)\otimes \mathrm{H}_{j+l-k}(P) \arrow[]{r}{\times} &
		\mathrm{H}_m(P\times P) \arrow[]{d}{\tau_C\cap} &
		\mathrm{H}_{m-k}(U) \arrow[]{d}{(c \mathrm{R})_*}
		\\
		& \mathrm{H}_{m-k}(U) \arrow[]{r}{(c \mathrm{R})_*} &
		\mathrm{H}_{m-k}(P) 
	\end{tikzcd}
	$$
	with $m = i+j+l-k$.
	We claim that this diagram commutes.
	The commutativity of the two squares which are not labelled is clear.
	The square labelled (1) commutes by the identity relating the cap and the cross-product, see \cite[A.3]{hingston:2017}.
	With the same identity we see that the square labelled (2) commutes up to the sign $(-1)^{ki}$.
	We are left with considering subdiagram (3) to prove this part of theorem.
	We shall use the naturality of the relative cap product, see \cite[Lemma A.1]{hingston:2017}.
	Let $X\in\mathrm{H}_{\bullet}(P^3)$ and consider
	\begin{equation}\label{eq_naturality_cap_ass_1}
		(c \mathrm{R})_* \big[  \tau_C \cap ((1\times c \mathrm{R})_* ((1\times \tau_C)\cap X ) )  \big]   =  (c \mathrm{R})_* (1\times c \mathrm{R})_* \big[     (1\times c \mathrm{R})^* \tau_C \cap ((1\times \tau_C)\cap X)       \big]  
	\end{equation}
	where we use the naturality property \cite[Lemma A.1]{hingston:2017}.
	Note that the pullback $(1\times c R)^*\tau_C$ is understood as a pullback along the map
	$$   (1\times c \mathrm{R})\colon (U_1,U_{1,\geq\epsilon_0}) \to  (U_{C^f},U_{C^f,\geq\epsilon_0})     $$
	with
	$$  (U_1,U_{1,\geq\epsilon_0}) =  \big((1\times c \mathrm{R})|_{P\times U}\big)^{-1}(U_{C^f},U_{C^f,\geq\epsilon_0})    $$
	Recall that $\tau_C = (p_1\times p_0)^*\tau_N$.
	One directly sees that 
	$$  (p_1\times p_0)\circ (1\times c \mathrm{R}) =    \mathrm{pr}_{1,2} \circ (p_1\times p_0)$$
	where $\mathrm{pr}_{1,2}\colon  P^3\to P^2 $ is the projection on the first two factors.
	Therefore we have
	$$   (1\times c \mathrm{R})^* \tau_C =  \tau_C \times 1 .     $$
	Now, consider the expression we obtain from going around the right hand side of subdiagram (3).
	We have again by naturality
	\begin{equation}\label{eq_naturality_cap_ass_2}
		(c \mathrm{R})_*\big[   \tau_C\cap (c \mathrm{R}\times 1)_* ((\tau_C\times 1)\cap X)      \big] =  (c \mathrm{R})_*( c \mathrm{R}\times 1)_*\big[     (c \mathrm{R}\times 1)^*\tau_C \cap ((\tau_C\times 1)\cap X)    \big] . 
	\end{equation}
	Similarly as before one can argue that $(c \mathrm{R}\times 1)^*\tau_C = 1\times \tau_C$.
	Note that one interprets $c \mathrm{R}\times 1$ as a map of pairs $(U_2,U_{2,\geq\epsilon_0})\to (U_{C^f},U_{C^f,\geq\epsilon_0})$ with
	$$     (U_2,U_{2,\geq\epsilon_0}) =  \big((c \mathrm{R}\times 1)|_{U\times P}\big)^{-1} (U_{C^f},U_{C^f,\geq\epsilon_0})  .          $$
	By explicitly writing down $U_1$ and $U_2$ we see that we have
	$$   U_1 = U_2 \quad \text{and}\quad U_{1,\geq\epsilon_0} = U_{2,\geq\epsilon_0} .   $$
	It can also be seen without problems that the maps
	$$    (c \mathrm{R})\circ (1\times c \mathrm{R}) \colon (U_1,U_{1,\geq\epsilon_0})\to P   \quad \text{and} \quad ( c \mathrm{R})\circ (c \mathrm{R}\times 1)  \colon (U_2,U_{2,\geq\epsilon_0})\to P  $$
	are homotopic.
	Thus by equations \eqref{eq_naturality_cap_ass_1} and \eqref{eq_naturality_cap_ass_2} as well as by the graded commutativity of the cup product we see that diagram (3) commutes up to the sign $(-1)^{k^2} = (-1)^k$.
	Reintroducing the signs of the path product one then sees that the above diagram does indeed show the associativity of the path product.
	
	We now show that the path product is unital.
	Let $[N]\in\mathrm{H}_k(N)$ be the fundamental class of $N$ and let $s\colon N\to P^f$ be the map assigning to a point $p\in N$ the point $(p,p,c_{f(p)})\in P^f$ where $c_q\in PM$ is the constant path at the point $q\in M$.
	Let $\overline{C}\subseteq N\times P^f$ be the space
	$$    \overline{C} = \{ (p,(x_0,x_1,\gamma))\in N\times P^f\,|\, p = x_0   \} .    $$
	Under the identification of $N$ with its image under $s$ we consider $\overline{C}$ as a subset of $C^f$ and we can restrict the tubular neighborhood and the Thom class $\tau_{C^f}$ of $C^f$ to $\overline{C}$.
	We obtain a tubular neighborhood $U_{\overline{C}}$ and a restricted Thom class $\tau_{\overline{C}}$.
	It is then clear that the following diagram commutes
	$$  
	\begin{tikzcd}
		\mathrm{H}_{i}(N \times P^f) \arrow[]{r}{(s\times \mathrm{id})_*} \arrow[]{d}{} 
        & \mathrm{H}_{i}(P^f \times P^f) \arrow[]{d}{} \\
		\mathrm{H}_i( U_{\overline{C}}, U_{\overline{C}}\setminus \overline{C}) \arrow[]{r}{(s\times \mathrm{id})_*} \arrow[swap]{d}{\tau_{\overline{C}}\cap} & \mathrm{H}_i(U_{C^f},U_{C^f}\setminus C^f) \arrow[]{d}{\tau_{C^f}\cap}
		\\
		\mathrm{H}_{i-k}(U_{\overline{C}}) \arrow[]{r}{(s\times \mathrm{id})_*} \arrow[swap]{d}{(c \circ \mathrm{R}_{C^f})_*} & \mathrm{H}_{i-k}(U_{C^f}) \arrow[]{d}{(c \circ \mathrm{R}_{C^f})_*}
		\\
		\mathrm{H}_{i-k}(P^f)\arrow[]{r}{=} & \mathrm{H}_{i-k}(P^f) .
	\end{tikzcd}
	$$
	Let $\varphi\colon \mathrm{H}_i(N\times P^f)\to \mathrm{H}_{i-k}(P^f)$ be the composition of the vertical maps on the left side of the above diagram.
	From the properties of the Thom isomorphism, one sees that we have
	$$    \varphi ([N]\times A) = A\quad \text{for all}\,\,\, A\in\mathrm{H}_{\bullet}(P^f) .      $$
	Thus, the above diagram shows that $s_*[N]\wedge A = A$ for all $A\in \mathrm{H}_{\bullet}(P^f)$ and this shows the unitality of the path product.

	For the second part of the theorem consider the following diagram.
	$$
	\begin{tikzcd}
		\mathrm{H}_i(P^f\times P^f) \arrow[]{d}{}\arrow[]{r}{(p\times p)_*} & [2.5em] \mathrm{H}_i(N\times N\times N\times N) \arrow[]{d}{}
		\\
		\mathrm{H}_i(P^f\times P^f,P^f\times P^f\setminus C^f) \arrow[]{r}{(p\times p)_*} \arrow[]{d}{\tau_{C^f}\cap} & \mathrm{H}_i(N^4, N^4 \setminus N\times \Delta N\times N) \arrow[]{d}{\pi_{2,3}^*\tau_N \cap}
		\\
		\mathrm{H}_{i-k}(U_{C^f}) \arrow[]{r}{(p\times p)_*} \arrow[]{d}{(c \circ \mathrm{R}_{C^f})_*}
		& 
		\mathrm{H}_{i-k}(N\times U_N\times N) \arrow[]{d}{(\mathrm{id}\times (h\circ \mathrm{R}_N)\times \mathrm{id})_*} 
		\\
		\mathrm{H}_{i-k}(P^f) \arrow[]{r}{ p_*}  & \mathrm{H}_{i-k}(N\times N) .
	\end{tikzcd}
	$$
	Here, $\pi_{2,3}\colon N^4 \to N^2$ is the projection on the second and third factor.
    We claim that this diagram commutes.
	The commutativity of the second square of this diagram can be seen by the naturality of the cap product and by noting that 
	$$    \pi_{2,3}\circ (p\times p) = p_1\times p_0  .$$
	The commutativity of the other two squares is clear.
	Now, let $X,Y\in\mathrm{H}_{\bullet}(P^f)$ and assume that $$ p_* X = a\times b \quad \text{and}\quad p_* Y = c\times d \quad \text{for} \,\,\, a,b,c,d\in \mathrm{H}_{\bullet}(N) . $$
	Then we have
	\begin{eqnarray*}
		p_* (X\wedge Y) & = & (-1)^{k-k|X|} p_* \big( ( c \circ \mathrm{R}_{C^f})_* (\tau_{C^f}\cap (X\times Y))\big) 
		\\
		& = & (-1)^{k-k|X|} (-1)^{k|a|} (\mathrm{id}\times (h\circ \mathrm{R}_N)\times \mathrm{id})_* (a\times (\tau_N \cap (b\times c))\times d) 
		\\ & = & 
		(-1)^{k-k|X|} (-1)^{k|a|} (-1)^{k-k|b|} \, a\times (b\overline{\cap} c) \times d 
		\\ & = & 
		\mu_{\beta}(p_* X \otimes p_* Y) .
	\end{eqnarray*}
	In the second equality we used the compatibility of the cap and the cross-product which gives the sign $(-1)^{k|a|}$.
	For the third and fourth equality recall the sign in the definition of the intersection product, see the discussion before Proposition \ref{prop_intersection_pr_from} and that $|a| + |b| = |X|$.
	This completes the proof.
\end{proof}

\begin{remark}
	Let $M$ be a closed manifold and let $N$ be  a closed oriented submanifold of $M$.
    Consider the path product on the homology of the space $P_N M$.
	\begin{enumerate}
		\item     Note that the second property of Theorem \ref{theorem_properties_path_product} is the analogue to the fact that the evaluation $\mathrm{ev}_{\Lambda}\colon \Lambda M\to M$ induces a ring morphism 
		$$   (\mathrm{ev}_{\Lambda})_* \colon \, (\mathrm{H}_{\bullet}(\Lambda M),\wedge_{\mathrm{CS}}) \xrightarrow[]{\hphantom{coni}} (\mathrm{H}_{\bullet}(M),\overline{\cap})    $$
        between the Chas-Sullivan ring and the intersection ring of a manifold $M$.
		Let $j\colon \Omega M\hookrightarrow \Lambda M$ be the inclusion of the fiber of the free loop fibration.
		For the Chas-Sullivan product it holds the Gysin map $j_!\colon \mathrm{H}_i(\Lambda M)\to \mathrm{H}_{i-n}(\Omega M)$ is an algebra morphism when $\mathrm{H}_{\bullet}(\Omega M)$ is equipped with the Pontryagin product, see \cite[Proposition 3.4]{chas:1999}.
		Let $i\colon \Omega M\hookrightarrow P_N M$ be the inclusion of the fiber of the fibration $p \colon P_N M\to N\times N$.
		This is the inclusion of a submanifold of codimension $2k$ where $k$ is the dimension of $N$.
		Therefore the Gysin map $i_!\colon \mathrm{H}_{l}(P_N M)\to \mathrm{H}_{l-2k}(\Omega M)$ cannot be a morphism of algebras for degree reasons.
		\item We also want to remark that, in contrast to the free loop fibration $\mathrm{ev}_{\Lambda}\colon \Lambda M\to M$, the fibration $p \colon P_N M\to N\times N$ does not have a section in general.
		In fact it is easy to see that there is a section if and only if the inclusion of the submanifold $i\colon N\hookrightarrow M$ is null-homotopic.
		Indeed, if $s\colon N\times N\to P_N M$ is a section of $p$, then we consider the restriction $$g = s|_{N\times \{p_0\}}\colon N\times \{p_0\}\to P_N M$$ for some point $p_0\in N$.
		Define a map $H\colon N\times I\to M$ given by
		$$   H(p,s) = g(p,p_0)(s) \quad \text{for}\,\,\,p\in N,\,\,s\in I .     $$
		One checks that this is a homotopy between the inclusion $i\colon N\hookrightarrow M$ and the constant map $p\mapsto p_0$.
		Conversely, if the inclusion $i\colon N\hookrightarrow M$ is null-homotopic then one can construct a section $N\times N\to P_N M$.
	\end{enumerate}
\end{remark}

\section{Examples of the path product}\label{sec_path_product_examples}

In this section we consider three situations where the path product can be described very explicitly.
We begin by studying the situation where $f\colon N\to M$ is the constant map.
One might wonder why we study the case that $f\colon N\to M$ is constant so extensively.
Once we have proven the invariance of the path product under homotopies of $f$ in the next section, this case will give us the behavior of the path product for all null-homotopic maps $f$.
Afterwards we will consider the case when $f\colon N\to M$ is a fiber bundle and see that this case actually gives an operation on the homology of a closed manifold.
Lastly, we study the diagonal map $\Delta\colon M\to M\times M$ for a closed oriented manifold $M$ and shall see that this recovers the Chas-Sullivan product on $M$.

\subsection{The constant map}
Let $M$ be a closed manifold and let $$\star\colon \mathrm{H}_i(\Omega M) \otimes \mathrm{H}_j(\Omega M)\to \mathrm{H}_{i+j}(\Omega M)$$
be the Pontryagin product, i.e.
$$  X\star Y = \mathrm{concat}_* \, (X\times Y) \quad \text{for} \quad  X,Y\in\mathrm{H}_{\bullet}(\Omega M) .    $$
Moreover, recall that we have the product $\mu_{\beta}$ induced by the intersection product form.
We define a product 
$$\mu_{N,\Omega}\colon \big(\mathrm{H}_{\bullet}(N)\otimes \mathrm{H}_{\bullet}(N)\otimes \mathrm{H}_{\bullet}(\Omega M)\big)^{\otimes2} \to  \mathrm{H}_{\bullet}(N)\otimes \mathrm{H}_{\bullet}(N)\otimes \mathrm{H}_{\bullet}(\Omega M) $$
by setting
$$   \mu_{N,\Omega} \big( (a\otimes b\otimes x) \otimes (c\otimes d\otimes y)\big) = \, (-1)^{|x|(|c|+|d|+k)} \mu_{\beta}(a\otimes b\otimes c\otimes d) \otimes (x\star y)     $$
for $a,b,c,d\in\mathrm{H}_{\bullet}(N)$ and $x,y\in\mathrm{H}_{\bullet}(\Omega M)$.
\begin{remark}
	Let $M$ be a closed manifold and $N$ a closed oriented manifold.
	Define an unsigned version of the product $\mu_{N,\Omega}$ by
	$$       \widetilde{\mu}\colon  \big(\mathrm{H}_{\bullet}(N)\otimes \mathrm{H}_{\bullet}(N)\otimes \mathrm{H}_{\bullet}(\Omega M)\big)^{\otimes2} \to  \mathrm{H}_{\bullet}(N)\otimes \mathrm{H}_{\bullet}(N)\otimes \mathrm{H}_{\bullet}(\Omega M)    $$
	with
	$$   \widetilde{\mu} \big( (a\otimes b\otimes x) \otimes (c\otimes d\otimes y)\big) = \, \mu_{\beta}(a\otimes b\otimes c\otimes d) \otimes (x\star y)  = \beta(b\otimes c) (a\otimes d \otimes x\star y)    $$
	for $a,b,c,d\in\mathrm{H}_{\bullet}(N)$ and $x,y\in\mathrm{H}_{\bullet}(\Omega M)$.
	Since the product $\mu_{\beta}$ and the Pontryagin product are associative one checks directly that $\widetilde{\mu}$ is associative as well.
	It might thus at first seem surprising that the pairing $\mu_{N,\Omega}$ with the additional sign is also associative.
	For 
	$$    a,b,c,d,e,f\in \mathrm{H}_{\bullet}(N) \quad \text{and}\quad x,y,z\in\mathrm{H}_{\bullet}(\Omega M)   $$
	one computes that
	$$   \mu_{N,\Omega}\big(\,a\otimes b\otimes x\otimes\, \mu_{N,\Omega}(c\otimes d\otimes y\otimes  e\otimes f\otimes z)\big)  $$  and $$  \mu_{N,\Omega}\big(\,\mu_{N,\Omega}(a\otimes b\otimes x\otimes c\otimes d\otimes y)\, \otimes  e\otimes f\otimes z)\big)     $$
	differ by a sign $(-1)^{|x|(|d|+|e|+k)}$.
	However, the pairing $\beta(d\otimes e)$ is zero if $|d|+|e| \neq k$, therefore this sign is indeed always $1$ for homology classes with non-trivial product $\mu_{N,\Omega}$.
	Consequently, the product $\mu_{N,\Omega}$ is associative.
\end{remark}

\begin{prop}\label{prop_trivial_map_path_space}
	Let $M$ be a closed manifold and $N$ a closed oriented manifold.
	Let $f\colon N\to M$ be the constant map $f(x) = p_0$ for the basepoint $p_0\in M$.
	Take homology with coefficients in a field.
	Then the homology $(\mathrm{H}_{\bullet}(P^f),\wedge)$ equipped with the path product is isomorphic as an algebra to $$ ( \mathrm{H}_{\bullet}(N)\otimes \mathrm{H}_{\bullet}(N)\otimes \mathrm{H}_{\bullet}(\Omega M), \mu_{N,\Omega} ) .  $$
\end{prop}
\begin{proof}
	Throughout this proof we shall write $\Omega$ for the based loop space $\Omega M$.
	Let $f\colon N\to M$ be the constant map to the basepoint of $M$.
	One checks directly that this yields $P^f = N\times N\times \Omega $.
	Moreover, note that $C^f$ is then the space
	$$     C^f  =  \{ ((x_0,x_1,\gamma) , (y_0,y_1,\sigma)  ) \in N^2\times \Omega \times N^2\times \Omega  \,|\,  x_1 = y_0  \}  \cong N^3 \times \Omega ^2.    $$
	Let $T\colon P^f\times P^f\to N^4\times \Omega^2$ be the swapping map
	$$    T((x_0,x_1,\gamma),(y_0,y_1,\sigma)) =  (x_0,x_1,y_0,y_1,\gamma,\sigma)      $$
	for 
	$(x_0,x_1,\gamma),(y_0,y_1,\sigma)\in P^f$.
	Note that $T$ maps the pair $(U_{C^f},U_{C^f}\setminus C^f)$ to the pair
	$$      (N\times U_N\times N \times \Omega^2 ,(N\times U_N\times N \times \Omega^2) \setminus( N\times \Delta N\times N\times \Omega^2) ) .      $$
	Moreover, let $\pi_{2,3}\colon N^4\times \Omega^2\to N^2$ be the projection on the second and third factor and let $h\colon \Delta N\to \{\mathrm{pt}\}$ be the trivial map.
		
	In the following diagram we denote a topological pair of the form $(X,X\setminus A)$ by $(X,\sim A)$.
	We claim that the following diagram commutes.
	$$
	\begin{tikzcd}
		\mathrm{H}_i(P^f\times P^f) \arrow[]{d}{} \arrow[]{r}{T_*}
		&
		\mathrm{H}_i(N^4 \times \Omega^2) \arrow[]{d}{}
		\\
		\mathrm{H}_i(U_{C^f}, \sim C^f) \arrow[]{d}{\tau_{C^f} \cap} \arrow[]{r}{T_*} 
		&
		\mathrm{H}_i(N\times U_N\times N\times \Omega^2, \sim N\times \Delta N\times N\times \Omega^2) 
		\arrow[]{d}{\pi_{2,3}^*\tau_N\cap}
		\\
		\mathrm{H}_{i-k}(U_{C^f}) \arrow[]{d}{(\mathrm{R}_{C^f})_*} \arrow[]{r}{T_*}
		&
		\mathrm{H}_i(N\times U_N\times N\times \Omega^2)
		\arrow[]{d}{(\mathrm{id}\times \mathrm{R}_N\times \mathrm{id})_*}
		\\
	 \mathrm{H}_{i-k}(C^f) \arrow[]{d}{c_*} \arrow[]{r}{T_*}
		&
		\mathrm{H}_{i-k}(N\times \Delta N\times N\times \Omega^2) \arrow[]{d}{(\mathrm{id}\times h\times \mathrm{id}\times\mathrm{concat})_*}
		\\
		\mathrm{H}_{i-k}(P^f )  \arrow[]{r}{=}
		&
		\mathrm{H}_{i-k}(N\times N\times \Omega )
	\end{tikzcd}   
	$$ 

The commutativity of the first square is clear.
For the second square one uses naturality of the cap product.
The third square commutes since the underlying diagram of maps commutes up to homotopy.
Indeed, let $((x_0,x_1,\gamma),(y_0,y_1,\sigma))\in U_{C^f}$.
Then we have
$$     (\mathrm{id}\times R_N \times \mathrm{id})\circ T ((x_0,x_1,\gamma),(y_0,y_1,\sigma)) =   (x_0,x_1,x_1,y_1,\gamma,\sigma)   $$
and 
$$   T\circ R_{C^f} ((x_0,x_1,\gamma),(y_0,y_1,\sigma)) =  (x_0,x_1,x_1,\gamma_1,\gamma,\mathrm{concat}(c_{p_0}, \sigma)) $$
where $c_{p_0}\in \Omega M$ is the constant loop at the basepoint.
By reparametrization we find a homotopy between the map $$\Omega M\to \Omega M, \sigma \mapsto \mathrm{concat}(c_{p_0},\sigma)  $$
and the identity on $\Omega M$.
Using this homotopy we see that 
$$     (\mathrm{id}\times R_N \times \mathrm{id})\circ T  \simeq  T\circ R_{C^f}.  $$
For the commutativity of the fourth square one checks that the underlying diagram of maps commutes strictly.
Up to the sign $(-1)^{k-k|X|}$ the left hand side of the above diagram gives the path product $X\wedge Y$, i.e. we get
$$   X\wedge Y =  (-1)^{k-k|X|}  (\mathrm{id}\times (h\circ \mathrm{R}_N)\times \mathrm{id}\times \mathrm{concat})_*
\big( \pi_{2,3}^*\tau_N \cap ( T_* (a\times b \times x\times c\times d\times  y)   )\big) .  $$
In order to understand the right hand side of the above diagram consider the diagram
$$   
\begin{tikzcd}
	\mathrm{H}_{\bullet}(N^2\times \Omega \times N^2\times \Omega) \arrow[]{d}{T_*} 
	&
	\mathrm{H}_{\bullet}(N)^{\otimes 2}\otimes \mathrm{H}_{\bullet}(\Omega)\otimes
	\mathrm{H}_{\bullet}(N)^{\otimes 2}\otimes \mathrm{H}_{\bullet}(\Omega)\arrow[swap]{l}{\times} \arrow[]{d}{\mathcal{T}}
	\\
	\mathrm{H}_{\bullet}(N^4\times \Omega^2) \arrow[]{d}{{\pi}_{2,3}^*\tau_N\cap}
	&
	\mathrm{H}_{\bullet}(N)^{\otimes 4}\otimes \mathrm{H}_{\bullet}(\Omega)^{\otimes 2}
	\arrow[swap]{l}{\times} \arrow[]{dd}{\mathrm{id}\otimes\beta\otimes\mathrm{id}\otimes \star}
	\\
	\mathrm{H}_{\bullet}(N\times U_N\times N\times \Omega^2) \arrow[]{d}{(\mathrm{id}\times (h\circ \mathrm{R}_N)\times \mathrm{id}\times \mathrm{concat})_*}
	&
	\\
	\mathrm{H}_{\bullet}(N\times N\times \Omega) 
	&
	\mathrm{H}_{\bullet}(N)^{\otimes 2}\otimes \mathrm{H}_{\bullet}(\Omega) \arrow[swap]{l}{\times}
\end{tikzcd}
$$
Here, $\mathcal{T}$ is the respective swapping of the tensor factors.
Let 
$$     a\otimes b\otimes x\otimes c\otimes d\otimes y\in  \mathrm{H}_{\bullet}(N)^{\otimes 2}\otimes \mathrm{H}_{\bullet}(\Omega)\otimes
\mathrm{H}_{\bullet}(N)^{\otimes 2}\otimes \mathrm{H}_{\bullet}(\Omega) .    $$
Then the upper part of this diagram commutes up to the sign $(-1)^{(|c|+|d|)|x|}$
and the lower square commutes up to the sign $(-1)^{k-k(|a|+|b|)}$.
Moreover, the composition down the right hand side of the diagram is the map $\mu_{N,\Omega}$ up to sign $(-1)^{(|c|+|d|+k)|x|}$.
We then find that
\begin{eqnarray*}
	&  \mu_{N,\Omega} ((a\otimes b\otimes x)\otimes (c\otimes d\otimes y)) = \\ &    (-1)^{k-k(|a|+|b|+|x|)}  (\mathrm{id}\times (h\circ \mathrm{R}_N)\times \mathrm{id}\times \mathrm{concat})_*\big( \pi_{2,3}^*\tau_N \cap ( T_*(a\times b\times c\times d\times x\times y))  \big) ,  
\end{eqnarray*}
so we see that this agrees with $X\wedge Y$, since $|X| = |a| + |b|+ |x|$.
This completes the proof.    
\end{proof}

We shall see some examples of the product $\mu_{N,\Omega}$ in the next section.
The relevance of the product $\mu_{N,\Omega}$ will become more clear once we have shown the invariance of the path product.


\subsection{Fiber bundles}

 We now study how the path product behaves in the case that $f\colon N\to M$ is a fiber bundle, where the base, the fiber and the total space are smooth manifolds.
    For the reader's convenience we shall write $E$ instead of $N$ in order to stay closer to the more common notation for fibrations.
    In the situation of $f\colon E\to M$ being a fiber bundle the path space $P^f$ is homotopy equivalent to the fiber product $E\times_M E$.
    This can be seen by adapting the argument that the fibrant replacement of a fibration is the fibration itself to our situation, see e.g. \cite[Section 5.7]{tomDieck:2008}.
    Since we shall work with the details of this construction later on we give a detailed proof of this fact in the next lemma.
     \begin{lemma}\label{lemma_hom_equivalence_fiber_bundles}
	  	Let $f\colon E\to M$ be a fibration and consider the associated space $P^f$.
	  	There is a fiber homotopy equivalence $\Phi\colon P^f \to E\times_M E$ with respect to the fibrations $p_1\colon P^f\to E$ and $\mathrm{pr}_2\colon E\times_M E\to E$.
        Similarly, there is a fiber homotopy equivalence $\Phi'\colon P^f\to E\times_M E$ with respect to the fibrations $p_0\colon P^f\to E$ and $\mathrm{pr}_1\colon E\times_M E\to E$.
	  \end{lemma}
	  \begin{proof}
	  	We consider the following diagram.
	  	$$
	  		\begin{tikzcd}
	  			P^f \times \{0\} \arrow[]{r}{p_0} \arrow[hook]{d}{} & [1.5em] E \arrow[]{d}{f} 
	  			\\
	  			P^f\times I \arrow[]{r}{h} & M .
	  		\end{tikzcd}
	  	$$
	  	The map $h\colon P^f\times I$ is defined by $  h(p,q,\gamma,s) = \gamma(s) $. 
	  	One verifies that the diagram commutes and since $f$ is a fibration we obtain a lift
	  	$  H \colon P^f\times I\to E $.
	  	In particular we have 
	  	$$   f(H(p,q,\gamma,1)) = h(p,q,\gamma,1) = \gamma(1) = f(q) \quad \text{for all}\,\,\, (p,q,\gamma)\in P^f.$$
	  	We define a map $\Phi\colon P^f\to E\times_M E$ by setting
	  	$$       \Phi(p,q,\gamma) = \big( H(p,q,\gamma,1),q    \big) \quad \text{for} \,\,\, (p,q,\gamma)\in P^f .         $$
	  	We claim that this is a homotopy equivalence.
	  	Consider the map $s\colon E\times_M E\to P^f$ defined by
	  	$$    s(p,q) = (p,q, c_{f(q)}) \quad \text{for} \,\,\, (p,q)\in E\times_M E .       $$
	  	Here, $c_r$ is the constant path in $M$ at $r\in M$.
	  	Note that since $(p,q)\in E\times_M E$ we have that $f(p)=f(q)$ and consequently the map $s$ is well-defined.
	  	Consider the composition $\Phi\circ s\colon E\times_M E\to E\times_M E$.
	  	We have
	  	$$     \Phi\circ s (p,q) = \big( H(p,q,c_{f(q)},1),q\big) \quad \text{for} \,\,\,(p,q)\in E\times_M E.      $$
	  	A homotopy between $\Phi\circ s$ and $\id_{E\times_M E}$ is given by
	  	$$   E\times_M E\times I\ni (p,q,u) \mapsto     \big(  H(p,q,c_{f(p)},u),q  \big) \in E\times_M E     $$ 
	  	since $H(p,q,c_{f(p)},u) \in f^{-1}(\{f(q)\})$ by the choice of $H$.
	  	One argues similarly for the composition $s\circ \Phi$.
      	Moreover, the fibrations $p_1\colon P^f\to E$ and $\mathrm{pr}_2\colon E\times_M E\to E$ are respected by the homotopies.
        For the second part we consider the diagram
        $$
	  		\begin{tikzcd}
	  			P^f \times \{1\} \arrow[]{r}{p_1} \arrow[hook]{d}{} & [1.5em] E \arrow[]{d}{f} 
	  			\\
	  			P^f\times I \arrow[]{r}{h'} & M .
	  		\end{tikzcd}
	  	$$
        where $h'(p,q,\gamma,s) = \gamma(s)$.
        By the homotopy lifting property we obtain a map $H'\colon P^f\times I\to E$ which fits into the above diagram.
        In particular we have $f(H'(p,q,\gamma,s)) = \gamma(s)$.
        We define a map $\Phi'\colon P^f\to E\times_M E$ by setting
        $$    \Phi'(p,q,\gamma) =   (p,H'(p,q,\gamma,0))  \quad \text{for all}\,\,\,(p,q,\gamma)\in P^f .   $$
        As before one checks that $\Phi'\circ s$ is homotopic to the identity on $E\times_M E$ and that $s\circ \Phi'$ is homotopic to the identity on $P^f$.
        Moreover, by construction, the fibrations $p_0\colon P^f\to E$ and $\mathrm{pr}_1\colon E\times_M E\to E$ are respected.
	  \end{proof}

    Under the homotopy equivalences of the previous Lemma the space of concatenable paths $C^f$ is taken to the fiber product
    $$  C_{{E\times_M E}} := (E\times_M E )
    \tensor[{_{\mathrm{pr}_2}}]{{\times}}{_{\mathrm{pr}_1}}
    (E\times_M E) =  \{ ((p,q),(r,s))\in (E\times_M E)^2 \,|\, q = r\}    .  $$
    This is a codimension $k$ submanifold of $(E\times_M E)^2$ where $k = \mathrm{dim}(E)$.
    Define a map $$\pi\colon E\times_M E\times E\times_M E  \to E\times E \quad \text{by}\quad \pi((p,q),(r,s)) = (q,r). $$
    We pull back the tubular neighborhood $U_E\subseteq E\times E$ to obtain an open set
    $ \widetilde{U}_E := \pi^{-1}(U_E)  .   $
    The Thom class $\tau_E\in \mathrm{H}^k(U_E,U_E\setminus \Delta E)$ can be pulled back along the map $\pi$ to obtain a class $\widetilde{\tau}_E = \pi^* \tau_E \in \mathrm{H}^k(\widetilde{U}_E , \widetilde{U}_E \setminus C_{E\times_M E})$. 
    There is further a retraction $\mathrm{R}_{E\times_M E}\colon \widetilde{U}_E\to C_{E\times_M E}$ given by $$\mathrm{R}_{E\times_M E}((p,q),(r,s)) = ((p,q),(q,s)) .$$
    We introduce the map
    $$    c_{E\times_M E}\colon C_{E\times_M E} \to E\times_M E  ,\quad ((p,q),(q,s))\mapsto (p,s) . $$
    Define a product $\wedge_E \colon \mathrm{H}_{\bullet}(E\times_M E)\otimes \mathrm{H}_{\bullet}(E\times_M E) \to \mathrm{H}_{\bullet}(E\times_M E)$ as the composition
    \begin{eqnarray*}
   \wedge_{E\times_M E} \colon          
   \mathrm{H}_i(E\times_M E)\otimes \mathrm{H}_j(E\times_M E) 
   & \xrightarrow[]{\hphantom{}(-1)^{k-ki}\times\hphantom{}}  
   &  \mathrm{H}_{i+ j}((E\times_M E)^2) \\
        & \xrightarrow[]{\hphantom{blaibalabla}} & \mathrm{H}_{i+j}((E\times_M E)^2, (E\times_M E)^2 \setminus C_{E\times_M E}) 
        \\
        & \xrightarrow[]{\hphantom{il}\text{excision}\hphantom{l}} & \mathrm{H}_{i+j}(\widetilde{U}_E,\widetilde{U}_E\setminus C_{E\times_M E})
        \\
            & \xrightarrow[]{\hphantom{blii}\widetilde{\tau}_E\cap\hphantom{bil}} & \mathrm{H}_{i+j-k}(\widetilde{U}_E)
        \\
        & \xrightarrow[]{(\mathrm{R}_{E\times_M E})_*} & \mathrm{H}_{i+j-k}(C_{E\times_M E})
            \\ & \xrightarrow[]{\hphantom{}(c_{E\times_M E})_*\hphantom{}} & \mathrm{H}_{i+j-k}(E\times_M E) .    
\end{eqnarray*}
    
  \begin{prop}
	  	Let $f\colon E\to M$ be a fiber bundle with $E$ and $M$ closed oriented manifolds.
	  	Then for $X,Y\in\mathrm{H}_{\bullet}(P^f)$ it holds that
	  	$$    X\wedge Y =  s_* \big( \Phi_* X \wedge_E \Phi'_* Y\big)     $$
	  	where the maps $s\colon E\times_M E\to P^f$ as well as $\Phi,\Phi'\colon P^f\to E\times_M E$ are the homotopy equivalences of Lemma \ref{lemma_hom_equivalence_fiber_bundles}.
	  \end{prop}
    \begin{proof}
        We begin by claiming that the following diagram commutes.
        $$
        \begin{tikzcd}
            \mathrm{H}_{i}(P^f) \otimes \mathrm{H}_{j}(P^f) \arrow[]{r}{\Phi_*\otimes \Phi'_*}
            \arrow[]{d}{}
            & [2.5em]
            \mathrm{H}_{i}(E\times_M E)\otimes \mathrm{H}_{j}(E\times_M E) \arrow[]{d}{}
            \\
            \mathrm{H}_{i+j}(U_{C^f}, U_{C^f}\setminus C^f) \arrow[]{r}{(\Phi\times \Phi')_*}
            \arrow[]{d}{\tau_{C^f}\cap} 
            &
            \mathrm{H}_{i+j}(\widetilde{U}_E, \widetilde{U}_E \setminus C_{E\times_M E}) 
            \arrow[]{d}{\widetilde{\tau}_E\cap}
            \\
            \mathrm{H}_{i+j-k}(U_{C^f}) \arrow[]{r}{(\Phi\times \Phi')_*}
            \arrow[]{d}{(\mathrm{R}_{C^f})_*}
            &
            \mathrm{H}_{i+j-k}(\widetilde{U}_E)\arrow[]{d}{( \mathrm{R}_{E\times_M E})_*}
            \\
            \mathrm{H}_{i+j-k}(C^f) \arrow[]{r}{(\Phi\times \Phi')_*}
            \arrow[]{d}{c_*} &
            \mathrm{H}_{i+j-k}(C_{E\times_M E}) \arrow[]{d}{(c_{E\times_M E})_*}
            \\
            \mathrm{H}_{i+j-k}(P^f) & \mathrm{H}_{i+j-k}(E\times_M E) \arrow[swap]{l}{s_*}
        \end{tikzcd}
        $$
        The commutativity of the first square is easy to see.
        The second square commutes by naturality of the cap product.
        Indeed, for $X\in \mathrm{H}_{\bullet}(U_{C^f},U_{C^f}\setminus C^f)$ we have that
        \begin{equation}\label{equation_naturality_once_more}
                \widetilde{\tau}_E\cap ((\Phi\times \Phi')_* X) =   (\Phi\times \Phi')_* \big(    (\Phi\times \Phi')^* \widetilde{\tau}_E \cap X  \big)   .  
        \end{equation}
        Since $\widetilde{\tau}_E = \pi^* \tau_E$ we have that
        $$    (\Phi\times \Phi')^* \widetilde{\tau}_E  =   (\pi\circ (\Phi\times \Phi'))^* \tau_E   =  (p_2\times p_1)^*\tau_E =  \tau_{C^f}   $$
        since the identity $\pi\circ (\Phi\times \Phi') = p_1\times p_0$ holds by Lemma \ref{lemma_hom_equivalence_fiber_bundles}.
        For the third square note that $\mathrm{R}_{C^f}\colon U_{C^f}\to C^f$ is homotopic to a retraction, see Section \ref{sec_def_path_product}.
        Therefore, if $\iota\colon C^f\hookrightarrow U_{C^f}$ is the inclusion we have $\iota\circ \mathrm{R}_{C^f} \simeq \mathrm{id}_{U_{C^f}}$.
        Moreover, one checks that the compositions 
        $$   \mathrm{R}_{E\times_M E}\circ (\Phi\times \Phi')\circ \iota\circ \mathrm{R}_{C^f} \quad \text{and}\quad (\Phi\times \Phi')\circ \mathrm{R}_{C^f}     $$
        agree strictly. 
        Since $\iota\circ \mathrm{R}_{C^f} \simeq \mathrm{id}_{U_{C^f}}$ it follows that
        $$   \mathrm{R}_{E\times_M E}\circ (\Phi\times \Phi')\simeq  (\Phi\times \Phi')\circ \mathrm{R}_{C^f}   .  $$
        Consequently, the third square commutes.
        For the last square, let $((p,q,\gamma),(q,r,\sigma))\in C^f$, then it holds that
        $$    s\circ c_{E\times_M E}\circ (\Phi\times \Phi') ((p,q,\gamma),(q,r,\sigma)) = \big( H(p,q,\gamma,1),H'(q,r,\sigma,0), c_{f(q)}\big) .      $$
        The map $H'' \colon C^f \times I\to P^f$ given by
        $$  H'' ((p,q,\gamma),(q,r,\sigma),u) =   \big( H(p,q,\gamma,1-u),H'(q,r,\sigma,u),  \mathrm{concat}( \gamma|_{[1-u,1]}, \sigma|_{[0,u]}) \big)    $$
        is a homotopy between the map $c$ and the composition $ s\circ c_{E\times_M E}\circ (\Phi\times \Phi') $.
        Hence the fourth square commutes.
        This shows the commutativity of the diagram and the identity for the path product follows.
    \end{proof}

    We can apply this proposition to the case that $f = \mathrm{id}_M$ for a closed manifold $M$.
    In this case we have that $E= M$ and $E\times_M E = M\times_M M \cong M$.
    Moreover, we have that $C_{M\times_M M} \cong \Delta M$ and $\widetilde{\tau}_M = \tau_M$ under these identifications.
    We see that in this case the product $\wedge_E$ becomes to the intersection product on $M$.
    \begin{cor}
        Let $M$ be a closed oriented manifold of dimension $k$ and let $f\colon M\to M$ be the identity map. 
        Then the path space $P^f$ is canonically identified with the free path space
        $  PM      $.
        The maps $\Phi\colon PM\to M$ and $\Phi'\colon PM\to M$ given by
        $    \Phi(\gamma ) = \gamma(1) \quad \text{and} \quad \Phi'(\gamma) = \gamma(0)   $
        are homotopic homotopy equivalences and for all $X,Y\in\mathrm{H}_{\bullet}(PM)$ we have 
        $$      X\wedge Y =   s_* (\Phi_* X \,\, \overline{\cap} \,\, \Phi'_* Y)       $$
        where $\overline{\cap}\colon \mathrm{H}_{\bullet}(M)^{\otimes 2}\to \mathrm{H}_{\bullet-k}(M)$ is the intersection product and where $s$ is a homotopy inverse to the maps $\Phi$ and $\Phi'$.
    \end{cor}

\subsection{The diagonal map}
    We end this section by studying the path product in the case of the diagonal map $\Delta\colon M\to M\times M$ for a closed oriented manifold $M$.
    We will see that the path product can be expressed by the Chas-Sullivan product.
    First, we review the definition of the Chas-Sullivan product.
    Let $M$ be a closed oriented manifold with $\mathrm{dim}(M) = n$ and denote by $\Lambda M$ its free loop space, i.e.
    $$   \Lambda M = \mathrm{ev}^{-1}(\Delta M)    $$
with $\mathrm{ev}\colon PM\to M\times M$ as in equation \eqref{eq_ev} and $\Delta M\subseteq M\times M$ the diagonal.
The restriction of the evaluation map thus becomes a map $$\mathrm{ev}_{\Lambda}\colon \Lambda M\to M    ,\quad \mathrm{ev}_{\Lambda}(\gamma) = \gamma(0) .  $$
Throughout this section we shall frequently write $\Lambda$ instead of $\Lambda M$.
Fix a Riemannian metric on $M$ and denote the distance function by $\mathrm{d}_M\colon M\times M\to \mathbb{R}$.
As before let $U_M$ be the space
$$    U_M = \{(p,q)\in M\times M\,|\, \mathrm{d}_M(p,q)<\epsilon'  \} .    $$
This is a tubular neighborhood of the diagonal $\Delta M\subseteq M\times M$ and there is thus a Thom class $\tau_M\in \mathrm{H}^n(U_M,U_M\setminus M)$.
Consider the space $\Lambda\times_M \Lambda$ which is defined as the pullback
$$  \Lambda\times_M\Lambda = (\mathrm{ev}_{\Lambda}\times \mathrm{ev}_{\Lambda})^{-1}(\Delta M) = \{(\gamma,\sigma)\in \Lambda\times \Lambda\,|\,\gamma(0) = \sigma(0)\}\,  .     $$
We call this space the  \textit{figure-eight space} since a pair $(\gamma,\sigma)\in \Lambda\times_M \Lambda$ can be understood as a map of the figure-eight $\mathbb{S}^1\vee\mathbb{S}^1\to M$.
The pullback of the tubular neighborhood $U_M$ gives us a tubular neighborhood of the figure-eight space
$$   U_{\mathrm{CS}} = (\mathrm{ev}_{\Lambda}\times\mathrm{ev}_{\Lambda})^{-1}(U_M)  = \{(\gamma,\sigma)\in\Lambda\times\Lambda\,|\,\mathrm{d}(\gamma(0),\sigma(0))<\epsilon'  \}  \,   , $$
see e.g. \cite[Proposition 2.2]{hingston:2017}.
We pull back the Thom class $\tau_M$ along the map $\mathrm{ev}_{\Lambda}\times \mathrm{ev}_{\Lambda}$ and get a class
$$    \tau_{\mathrm{CS}} = (\mathrm{ev}_{\Lambda}\times\mathrm{ev}_{\Lambda})^* \tau_M \in \mathrm{H}^n(U_{\mathrm{CS}},U_{\mathrm{CS}}\setminus \Lambda\times_M \Lambda) \, .      $$
Let $R_{\mathrm{CS}}\colon U_{\mathrm{CS}}\to \Lambda\times_M\Lambda$ be the induced retraction of the tubular neighborhood. 
For an explicit description of $\mathrm{R}_{\mathrm{CS}}$ we follow \cite[Section 1]{hingston:2017} and we have for $(\gamma,\sigma)\in U_{\mathrm{CS}}$ that
$$\mathrm{R}_{\mathrm{CS}}(\gamma,\sigma) =  (\gamma, \mathrm{concat}(\widehat{\gamma(0)\sigma(0)}, \sigma, \widehat{\sigma(0)\gamma(0)})) . $$
Note that the concatenation of paths restricts to a map $\Lambda\times_M \Lambda\to \Lambda$.
We define the \textit{Chas-Sullivan product} as the composition 
\begin{eqnarray*}
	\wedge_{\mathrm{CS}} \colon  \mathrm{H}_{i}(\Lambda)\otimes \mathrm{H}_j(\Lambda) &\xrightarrow{  (-1)^{n-ni}\times   }& \mathrm{H}_{i+j}(\Lambda\times \Lambda)
	\\ &\xrightarrow{ \hphantom{coninci}\vphantom{r}\hphantom{conci} }& \mathrm{H}_{i+j}(\Lambda\times \Lambda, \Lambda\times \Lambda\setminus \Lambda\times_M \Lambda) \\
	&\xrightarrow{\hphantom{cii}\text{excision}\hphantom{ic}}&
	\mathrm{H}_{i+j}(U_{\mathrm{CS}}, U_{\mathrm{CS}}\setminus \Lambda\times_M\Lambda) \\
	&\xrightarrow{\hphantom{coci} \tau_{\mathrm{CS}}\cap \hphantom{coc} }& \mathrm{H}_{i+j-n}(U_{\mathrm{CS}})
	\\ &\xrightarrow[]{\hphantom{cin}(\mathrm{R}_{\mathrm{CS}})_* \hphantom{co}} & \mathrm{H}_{i+j-n}(\Lambda\times_M \Lambda) 
	\\ & \xrightarrow[]{\hphantom{ici}\mathrm{concat}_*\hphantom{ici}} & \mathrm{H}_{i+j-n}(\Lambda) \, .
\end{eqnarray*}

Now, consider the diagonal map $\Delta\colon M\to M\times M$ for a closed $n$-dimensional manifold $M$.
Note that for an element $(p,q,\gamma)\in P^{\Delta}$ we have that $\gamma = (\gamma_a,\gamma_b)\in P(M\times M)$ such that $\gamma_a(0) = \gamma_b(0) = p$ and $\gamma_a(1) =\gamma_b(1) = q$.
We define a map $\varphi\colon P^{\Delta}\to \Lambda M$ by 
$$    \varphi (p,q,(\gamma_a,\gamma_b))  =  \mathrm{concat}(\gamma_a,\overline{\gamma_b}) .   $$
where for a path $\eta\in PM$ the notation $\overline{\eta}$ means the reversed path, i.e. $\overline{\eta}(t) = \eta(1-t)$.
Similarly we define $\varphi'\colon P^{\Delta}\to \Lambda$ by setting
$$    \varphi'(p,q,(\gamma_a,\gamma_b)) =   \mathrm{concat}(\overline{\gamma_b},\gamma_a) . $$

\begin{lemma}\label{lemma_homeomorphisms_free_loop_space}
    The maps $\varphi\colon P^{\Delta}\to \Lambda M$ and $\varphi'\colon P^{\Delta}\to \Lambda M$ are homeomorphisms.
    Furthermore, the identities $p_0 = \mathrm{ev}_{\Lambda}\circ \varphi $ and $p_1 = \mathrm{ev}_{\Lambda}\circ \varphi'$ hold and the homeomorphisms $\varphi$ and $\varphi'$ are homotopic to each other.
\end{lemma}
\begin{proof}
   
    One defines an inverse to $\varphi$ by defining
    $$  \chi\colon \Lambda M\to P^{\Delta},\quad \chi(\gamma) = (\gamma(0),\gamma(\tfrac{1}{2}), (\gamma|_{[0,\tfrac{1}{2}]}, \overline{\gamma|_{[\tfrac{1}{2},1]}} )) .  $$ 
    One checks that this indeed an inverse to $\varphi$ and that the identity $p_0 = \mathrm{ev}_{\Lambda}\circ \varphi$ holds.
    Similarly one verifies that the map
    $$   \chi'\colon \Lambda M\to P^{\Delta}, \quad \chi'(\gamma) =  (\gamma(\tfrac{1}{2}),\gamma(0), (\gamma|_{[\tfrac{1}{2},1]}, \overline{\gamma|_{[0,\tfrac{1}{2}]}}))     $$
    is an inverse to $\varphi'$ and that $p_1 = \mathrm{ev}_{\Lambda}\circ \varphi'$.
    The homotopy between $\varphi$ and $\varphi'$ is obtained by noting that the map $\beta_{\!\frac{1}{2}}\colon \Lambda M\to \Lambda M$, $\gamma\mapsto \gamma(t+ \tfrac{1}{2})$ is homotopic to the identity and that $\varphi$ can be expressed as the composition $\varphi = \beta_{\!\frac{1}{2}}\circ \varphi'$.
\end{proof}
\begin{prop}\label{prop_path_product_diagonal_map}
    Let $M$ be a closed oriented manifold and let $\Delta\colon M\to M\times M$ be the diagonal map.
    Then the path product for $P^{\Delta}$ can be expressed by the Chas-Sullivan product.
    More precisely, for $X,Y\in\mathrm{H}_{\bullet}(P^{\Delta})$ we have that
    $$       \varphi_* (X\wedge Y )  =   \varphi_* X \wedge_{\mathrm{CS}} \varphi_* Y  .        $$
\end{prop}
\begin{proof}
    We show that the following diagram commutes.
     $$
        \begin{tikzcd}
            \mathrm{H}_{i}(P^{\Delta}) \otimes \mathrm{H}_{j}(P^{\Delta}) \arrow[]{r}{\varphi'_*\otimes \varphi_*}
            \arrow[]{d}{}
            &
            \mathrm{H}_{i}(\Lambda M)\otimes \mathrm{H}_{j}(\Lambda M) \arrow[]{d}{} 
            \\
            \mathrm{H}_{i+j}(U_{C^{\Delta}}, U_{C^{\Delta}} \setminus C^{\Delta}) \arrow[]{r}{(\varphi'\times \varphi)_*}
            \arrow[]{d}{\tau_{C^{\Delta}}\cap}
            & 
            \mathrm{H}_{i+j}(U_{\mathrm{CS}}, U_{\mathrm{CS}}\setminus \Lambda \times_M \Lambda)\arrow[]{d}{\tau_{\mathrm{CS}}\cap}
            \\
            \mathrm{H}_{i+j-n}(U_{C^{\Delta}}) \arrow[]{r}{(\varphi'\times \varphi)_*}
            \arrow[]{d}{(\mathrm{R}_{C^{\Delta}})_*}
            &
            \mathrm{H}_{i+j-n}(U_{\mathrm{CS}})\arrow[]{d}{( \mathrm{R}_{\mathrm{CS}})_*}
            \\
            \mathrm{H}_{i+j-n}(C^{\Delta}) \arrow[]{r}{(\varphi'\times \varphi)_*}
            \arrow[]{d}{c_*} &
            \mathrm{H}_{i+j-n}(\Lambda\times_M \Lambda) \arrow[]{d}{\mathrm{concat}_*}
            \\
            \mathrm{H}_{i+j-n}(P^{\Delta}) \arrow[]{r}{\varphi'_*} & \mathrm{H}_{i+j-n}(\Lambda M) .
        \end{tikzcd}
        $$
        The commutativity of the first square is easy to see.
        For the second square, let $X\in\mathrm{H}_{\bullet}(U_{C^{\Delta}}, U_{C^{\Delta}}\setminus C^{\Delta})$.
        Then we have by naturality that
        $$  \tau_{\mathrm{CS}}\cap (\varphi'\times \varphi)_* X =   (\varphi'\times \varphi)_* \big( (\varphi'\times \varphi)^* \tau_{\mathrm{CS}} \cap X\big)  .    $$
        Recall that $\tau_{\mathrm{CS}} = (\mathrm{ev}_{\Lambda}\times \mathrm{ev}_{\Lambda})^* \tau_M$.
        Hence,
        $$    (\varphi'\times \varphi)^* \tau_{\mathrm{CS}} =   ((\mathrm{ev}_{\Lambda}\times \mathrm{ev}_{\Lambda})\circ (\varphi'\times\varphi))^* \tau_M =  (p_1\times p_0)^*\tau_M = \tau_{C^{\Delta}}$$
        where we used Lemma \ref{lemma_homeomorphisms_free_loop_space} for the identity 
        $$    (\mathrm{ev}_{\Lambda}\times \mathrm{ev}_{\Lambda}) \circ (\varphi\times \varphi') = p_1\times p_0 .    $$
        Consequently, the second square commutes.
        The third square commutes because the underlying diagram of maps commutes strictly.
        For the fourth square, we claim that the underlying diagram of maps commutes up to homotopy.
        Indeed, let $((p,q,\gamma),(q,r,\sigma))\in C^{\Delta}$, then we have
        $$  \varphi'\circ c((p,q,\gamma),(q,r,\sigma)) =  (p,r,\mathrm{concat}(\gamma_a,\sigma_a,\overline{\sigma_b}, \overline{\gamma_b}))     $$
        where we write $\gamma = (\gamma_a,\gamma_b)$ and $\sigma = (\sigma_a,\sigma_b)$. 
        On the other hand we have
        $$   \mathrm{concat}\circ (\varphi'\times \varphi) ((p,q,\gamma),(q,r,\sigma))  = \mathrm{concat}(\overline{\gamma_b},\gamma_a,\sigma_a,\overline{\sigma_b}) .  $$
        Define the map $\beta_{\!\frac{1}{4}}\colon \Lambda M\to \Lambda M$ by $\beta_{\!\frac{1}{4}}(\gamma)(t) =  \gamma(t-\tfrac{1}{4})$.
        Using the above expressions one checks that $$\beta_{\tfrac{1}{4}}\circ \varphi'\circ c = \mathrm{concat}\circ (\varphi'\times\varphi)  .  $$ 
        Since $\beta_{\!\frac{1}{4}}$ is homotopic to $\mathrm{id}_{\Lambda M}$ we conclude that $ \varphi'\circ c \simeq \mathrm{concat}\circ (\varphi'\times\varphi)$.
        This shows that the fourth square commutes.
        Recall from Lemma \ref{lemma_homeomorphisms_free_loop_space} that there is a homotopy $\varphi\simeq \varphi'$.
        By noting that up the sign $(-1)^{n-ni}$ the vertical compostions of the diagram at the beginning of the proof are the path product, respectively the Chas-Sullivan product, the claim follows.
\end{proof}

\section{Invariance of the path product}\label{sec_invariance_product}

In this section we show that if $f\colon N\to M$ and $g\colon N\to M$ are homotopic maps, then the respective path products on $\mathrm{H}_{\bullet}(P^f)$ and $\mathrm{H}_{\bullet}(P^g)$ are isomorphic.
We also discuss consequences of this result.

Let $f,g\colon N\to M$ be smooth maps between closed manifolds $N$ and $M$.
Assume that $H\colon N\times I\to M$ is a homotopy between $f$ and $g$, i.e.
$$    H(p,0) = f(p) \quad \text{and} \quad H(p,1) = g(p) \quad \text{for all}\,\,\,p\in N.     $$
Define $\eta\colon N\to PM$ by setting
$$     \eta(p)(s)  =  H(p,s) \quad \text{for all}\,\,\, p\in N,\,\, s\in I.     $$
We define a map $\Phi\colon  P^f \to P^g$ by
$$     \Phi(x_0,x_1,\gamma) =  (x_0,x_1,\mathrm{concat} ( \overline{\eta(x_0)}, \gamma, \eta(x_1) ))       \quad \text{for}\,\,\ (x_0,x_1,\gamma)\in P^f.  $$
We also define $\Psi\colon P^g \to P^f$ by 
$$      \Psi(x_0,x_1,\gamma) =  (x_0,x_1,\mathrm{concat} (\eta(x_0),\gamma,\overline{\eta(x_1)}))   \text{for}\,\,\ (x_0,x_1,\gamma)\in P^g .   $$
One checks that $\Phi$ and $\Psi$ are homotopy inverses to each other and thus the path spaces $P^f$ and $P^g$ are homotopy equivalent.
In particular this shows directly that if $f\colon N\hookrightarrow M$ is a submanifold such that $f$ is null-homotopic, then the space of paths in $M$ with endpoints in $N$ is homotopy equivalent to
$$    P_N M \simeq N\times N\times \Omega M .     $$
In the following theorem we shall see that the path product is invariant under the homotopy equivalence $\Phi\colon P^f\xrightarrow{\simeq} P^g $.
\begin{theorem}\label{theorem_invariance_path_product}
	Let $M$ be a closed manifold and let $N$ be a closed oriented manifold.
	Take homology with coefficients in a commutative ring $R$.
	Let $f\colon N\to M$ and $g\colon N\to M$ be homotopic maps and let $\Phi\colon P^f\to P^g$ be the homotopy equivalence as above.
	Then the induced map in homology
	$$     \Phi_*  \colon \mathrm{H}_{\bullet}(P^f) \to \mathrm{H}_{\bullet}(P^g)      $$
	is an isomorphism of rings.
\end{theorem}
\begin{proof}
	The claim will follow from the commutativity of the following diagram.
	$$
		\begin{tikzcd}
				\mathrm{H}_i(P^f\times P^f)   \arrow[]{r}{(\Phi\times \Phi)_*}
				\arrow[]{d}{}
				& [3em]
				\mathrm{H}_i(P^g\times P^g) \arrow[]{d}{}
				\\
				\mathrm{H}_i(U_{C^f}, U_{C^f}\setminus C^f) \arrow[]{r}{(\Phi\times \Phi)_*} \arrow[swap]{d}{\tau_{C^f}\cap }
				&
				\mathrm{H}_i(U_{C^g}, U_{C^g}\setminus C^g) \arrow[]{d}{\tau_{C^g}\cap }
				\\
				\mathrm{H}_{i-k}(U_{C^f}) \arrow[]{r}{(\Phi\times \Phi)_*}
				 \arrow[swap]{d}{(c\circ \mathrm{R}_{C^f})_*}
				&
				\mathrm{H}_{i-k}(U_{C^g}) \arrow[]{d}{(c\circ \mathrm{R}_{C^g})_*}
				\\
				\mathrm{H}_{i-k}(P^f ) \arrow[]{r}{\Phi_*} & \mathrm{H}_{i-k}(P^g)   .
		\end{tikzcd}
	$$
    The commutativity of the first square is clear.
    For the second square we need to show that
    $$    \tau_{C^g}\cap ((\Phi\times \Phi)_* X) =  (\Phi\times \Phi)_* (\tau_{C^f}\cap X)        $$
    for $X\in \mathrm{H}_{\bullet}(U_{C^f},U_{C^f}\setminus C^f)$.
    By naturality the left hand side of the above equation equals
    $$     \tau_{C^g}\cap ((\Phi\times \Phi)_* X)  =  (\Phi\times \Phi)_* \big( ((\Phi\times \Phi)^* \tau_{C^g}) \cap X \big) .   $$
    Recall that $\tau_{C^f} = (p_1\times p_0)^*\tau_N$ where $p_0,p_1\colon P^f\to N$ are the fibrations as described in Section \ref{sec_def_path_product}.
    Similarly, $\tau_{C^g} = (\overline{p}_1\times \overline{p}_0)^*\tau_N$ where $\overline{p}_0,\overline{p}_1\colon P^g\to N$ are denoted with a bar so that they can be distinguished from $p_0$ and $p_1$.
    One checks that the identity  
    $$  (\overline{p}_1\times \overline{p}_0)\circ (\Phi\times \Phi) = (p_1\times p_0)   $$
    holds and thus we see that
    $$   (\Phi\times \Phi)^* \tau_{C^g} =  \tau_{C^f}  . $$
    This shows the commutativity of the second square.
    For the third square we claim that the underlying diagram of maps commutes up to homotopy.
    By the definitions we see that the map $c\circ \mathrm{R}_{C^g}\circ (\Phi\times \Phi)$ is homotopic to a map $U_{C^f}\to P^g$ which maps
    \begin{equation}\label{eq_homotopic_comp_1}
          ((x_0,x_1,\gamma),(y_0,y_1,\sigma)) \mapsto   (x_0,y_1,\mathrm{concat}(  \overline{\eta(x_0)}, \gamma,\eta(x_1), g(\widehat{x_1 y_0}) , \overline{\eta(y_0)} , \sigma, \eta(y_1)   ))     
    \end{equation}
    for $((x_0,x_1,\gamma),(y_0,y_1,\sigma))\in U_{C^f}$.
    On the other hand the composition $\Phi\circ c\circ \mathrm{R}_{C^f}$ is homotopic to a map $U_{C^f}\to P^g$ given by
    \begin{equation}\label{eq_homotopic_comp_2}
          ((x_0,x_1,\gamma),(y_0,y_1,\sigma)) \mapsto   (x_0,y_1,\mathrm{concat}(  \overline{\eta(x_0)}, \gamma,f(\widehat{x_1 y_0}), \sigma, \eta(y_1)   ))       . 
    \end{equation}
    In order to see that the two maps described in equations \eqref{eq_homotopic_comp_1} and \eqref{eq_homotopic_comp_2} are homotopic, we consider the maps
    $$  \rho_1\colon  U_N \to PM , \quad (x_1,y_0) \mapsto  f(\widehat{x_1 y_0})     $$
    and $$    \rho_2\colon  U_N \to PM  \quad (x_1,y_0) \mapsto \mathrm{concat}(\eta(x_1),g (\widehat{x_1y_0}), \overline{\eta(y_0)}) . $$
    Note that for $(x_1,y_0)\in U_N$ both $\rho_1(x_1,y_0)$ as well as $\rho_2(x_1,y_0)$ are paths from $f(x_1)$ to $f(y_0)$.
    These two maps are homotopic via a map which keeps the endpoints of the paths fixed.
    Indeed, consider the homotopy
    $$   U_N\times I \to PM ,\quad  ((x_1,y_0),s) \mapsto \mathrm{concat}( \eta(x_1)|_{[0,1-s]}, H(\widehat{x_1 y_0},1-s), \overline{\eta(y_0)}|_{[s,1]} ) .    $$
    This is a homotopy between the map $\rho_2$ and the map
    $$   U_N\to P_M   ,\quad (x_1,y_0)\mapsto \mathrm{concat}( c_{f(x_1)}, f(\widehat{x_1 y_0}), c_{f(y_0)})      $$
    where $c_q$, $q\in M$ denotes the constant path at the point $q$.
    But clearly the latter map is homotopic to $\rho_1$ and thus we obtain a homotopy between $\rho_1$ and $\rho_2$.
    Using this homotopy and the expressions in equations \eqref{eq_homotopic_comp_1} and \eqref{eq_homotopic_comp_2} one can see that the compositions 
    $$ c\circ \mathrm{R}_{C^g}\circ (\Phi\times \Phi) \quad \text{and}\quad  \Phi\circ c\circ \mathrm{R}_{C^f}  $$ are homotopic.
    This completes the proof.
\end{proof}
\begin{cor}
	Let $N$ be a closed oriented manifold and let $i_1,i_2\colon N\to M$ be two homotopic embeddings of $N$ into a closed manifold $M$.
	Then the path product rings 
	$$       \big( \mathrm{H}_{\bullet}(P_N^{i_1} M), \wedge\big)  \quad \text{and} \quad       \big( \mathrm{H}_{\bullet}(P_N^{i_2} M), \wedge\big)     $$
	are isomorphic.
\end{cor}
By Proposition \ref{prop_trivial_map_path_space} and Theorem \ref{theorem_invariance_path_product} we also obtain the following corollary.
\begin{cor}\label{cor_nullhomotopic_product}
	Let $M$ be a closed manifold and $N$ be a closed oriented submanifold.
	Assume that the inclusion $i\colon N\hookrightarrow M$ is null-homotopic.
	Let $R$ be a commutative ring and take homology with coefficients in $R$.
	\begin{enumerate}
		\item There is an isomorphism of rings
		$$    \big(    \mathrm{H}_{\bullet}(P_N M) ,\wedge     \big)    \cong \big( \mathrm{H}_{\bullet}(P^{\mathrm{cst}}), \wedge \big)  $$
		where $P^{\mathrm{const}}  = N\times N\times \Omega M$ is the path space with respect to the constant map.
		\item If the coefficient ring $R$ is a field then there is an isomorphism
		$$     \big(    \mathrm{H}_{\bullet}(P_N M) ,\wedge     \big)    \cong \big(    \mathrm{H}_{\bullet}(N)\otimes \mathrm{H}_{\bullet}(N) \otimes \mathrm{H}_{\bullet}(\Omega M),\mu_{N,\Omega}       \big)      $$
		with $\mu_{N,\Omega}$ the product as in Proposition \ref{prop_trivial_map_path_space}.
	\end{enumerate}
\end{cor}
\begin{example}
    Corollary \ref{cor_nullhomotopic_product} yields a very explicit way of computing the path product in many particular cases.
    \begin{enumerate}
        \item Let $M = \mathbb{S}^n$ be a sphere.
         Every closed submanifold $N\hookrightarrow\mathbb{S}^n$ has null-homotopic inclusion, thus in the case of the manifold $M$ being a sphere we can compute the path product from the knowledge of the intersection ring of $N$ and the well-known Pontryagin ring of the based loop space of the sphere.
         \item Let $G$ be a compact Lie group and $N = T^r\hookrightarrow G$ be a maximal torus.
         If $G$ is simply connected, then the inclusion of the maximal torus is null-homotopic.
         Moreover the intersection ring of the torus is well-known, as is the Pontryagin ring of the based loop space for many compact Lie groups, see e.g. \cite{bott:1958b}.
    \end{enumerate}
\end{example}

\section{Module structure over the Chas-Sullivan ring}\label{sec_definition_module}

We now turn to the module structure of $\mathrm{H}_{\bullet}(P^f)$ over the Chas-Sullivan ring.
We refer to Section \ref{sec_path_product_examples} for the definition of the Chas-Sullivan product.
Let $M$ and $N$ be closed oriented manifolds and let $f\colon N\to M$ be a smooth map.
Define $D^f$ to be the space
$$    D^f = \{   (\gamma, (x_0,x_1,\sigma)) \in \Lambda M\times P^f \,|\, \gamma(0) = \sigma(0) = f(x_0)     \} .      $$
This space can be obtained as a pullback which is expressed by the following diagram.
$$
	\begin{tikzcd}
		D^f \arrow[]{r}{} \arrow[]{d}{} &  [2.5em] \Lambda M\times P^f \arrow[]{d}{\mathrm{ev}_{\Lambda}\times \mathrm{ev}_0 }
		\\
		\Delta M \arrow[]{r}{}  & M\times M.
	\end{tikzcd}
$$
Here, $\mathrm{ev}_0 \colon P^f\to M$ is the map
$$  \mathrm{ev}_0((x_0,x_1,\gamma)) = \gamma(0)  = f(x_0) .     $$
In general, $\mathrm{ev}_0$ is not a fibration.
Note that $\mathrm{ev}_{\Lambda}\colon \Lambda M\to M$ is a submersion, so it follows that $\mathrm{ev}_{\Lambda}\times \mathrm{ev}_0\colon \Lambda M\times P^f$ is transverse to the diagonal $\Delta M$.
Consequently, $D^f$ is a submanifold of $\Lambda M\times P^f$ of codimension $n$.
Moreover, the normal bundle of $D^f$ is
$$   E^f  =  \big(  \mathrm{ev}_{\Lambda}\times \mathrm{ev}_0|_{D^f}  \big)^* TM .     $$
One verifies that we have the identity
$$  \mathrm{ev}_{\Lambda}\times \mathrm{ev}_0|_{D^f} =  (f\circ p_0\circ \mathrm{pr}_2 )|_{D^f}    $$
where $\mathrm{pr}_2\colon \Lambda M \times P^f\to P^f$ is the projection on the second factor.
In particular we thus see that
$$  E^f  \cong (p_0\circ \mathrm{pr}_2)^* \, f^* TM  =  \{ (\gamma,(x_0,x_1,\sigma),v)\in D^f \times TM \,|\, v\in T_{f(x_0)} M\} .    $$
It will be crucial in the next section that $E$ is the pullback of the bundle $f^* TM \to N$.
Next, we describe a tubular neighborhood of $D^f$ in $\Lambda M\times P^f$.
Set
$$  U_{D^f} =   \{(\gamma,(x_0,x_1,\sigma))\in \Lambda M \times P^f \,|\,  \mathrm{d}_M ( \gamma(0),\sigma(0)) < \epsilon' \}    .  $$
The map $\mathrm{ev}_{\Lambda}\times \mathrm{ev}_0\colon \Lambda M\times P^f$ restricts to a map of pairs
$$      (U_{D^f}, U_{D^f}\setminus D^f ) \to (U_M,U_{M}\setminus \Delta M) .       $$
Consequently, we can pull back the Thom class $\tau_M$ to a class
$$  \tau_{D^f}  =   (\mathrm{ev}_{\Lambda}\times \mathrm{ev}_0)^* \tau_M \in \mathrm{H}^n(U_{D^f}, U_{D^f}\setminus D^f ) .   $$
Define a homotopy retraction $\mathrm{R}_{D^f}\colon U_{D^f}\to D^f$ by setting
$$    \mathrm{R}_{D^f} (\gamma,(x_0,x_1,\sigma)) = \big( 
  \mathrm{concat}(\widehat{\sigma(0)\gamma(0)}, \gamma,\widehat{\gamma(0)\sigma(0)}, (x_0,x_1,\sigma)  \big)        $$
where $\widehat{pq}\in PM$ is the unique length-minimizing geodesic between $p,q\in M$ with $\mathrm{d}_M(p,q)<\epsilon'$.
In the following lemma we prove that $U_{D^f}$ is indeed a tubular neighborhood for $D^f$.

\begin{lemma}\label{lemma_retraction_explicit}
    Let $M$ and $N$ be closed manifolds with a smooth map $f\colon N\to M$.
    \begin{enumerate}
        \item     The subset $U_{D^f}\subseteq \Lambda M\times P^f$ is a tubular neighborhood of $D^f$.
    More precisely, there is a homeomorphism $\nu_D\colon E^f \to U_{D^f}$ which restricts to the identity of $D^f$ on the zero-section of the bundle $E^f\to D^f$.
    \item The tubular neighborhood $U_{D^f}$ is compatible with the tubular neighborhood $U_M$, i.e. the following diagram commutes
    $$  \begin{tikzcd}
         E^f \arrow[]{r}{\nu_{D}} \arrow[]{d}{}  & [2em] \Lambda M\times P^f\arrow[]{d}{\mathrm{ev}_{\Lambda}\times \mathrm{ev}_0} \\ TM \arrow[]{r}{\nu_M} & M\times M 
    \end{tikzcd}     $$
    with $\nu_M\colon TM\to M\times M$ the tubular neighborhood given by $\nu_M(p,v) = (\exp_p(v),p)$ where $TM$ is identified with the open disk bundle 
    $  \{ v\in TM\,|\, |v| < \epsilon'   \}     $.
    \item The retraction $k_D \colon U_{D^f} \to D^f$ which is induced by the tubular neighborhood $\nu_{D}\colon E^f\xrightarrow[]{\cong} U_{D^f}$ is homotopic to the map $\mathrm{R}_{D^f}$.
    \end{enumerate}
\end{lemma}
\begin{proof}
     We closely follow the proofs of \cite[Proposition 2.2]{hingston:2017} and \cite[Lemma 2.3]{hingston:2017}.
    Recall that we have a chosen a Riemannian metric on $M$.
    We identify $TM$ with the open disk bundle
    $  \{ v\in TM\,|\, |v| < \epsilon'   \}     $.
    Note that a retraction $\mathrm{R}_M\colon U_M \to M$ is given by $\mathrm{R}_M(p,q) = q$ for $(p,q)\in U_M$.
     
    Recall that we denote the injectivity radius of $M$ by $\rho$.
    Hingston and Wahl define in \cite[Section 2]{hingston:2017} a map $h\colon U_M\times M\to M$ with the following properties:
    \begin{itemize}
    \item For $(u,v)\in U_M$ the map $h(u,v) = h(u,v,\cdot)\colon M\to M$ is a diffeomorphism.
    \item It holds that $h(u,v)(u) = v$ for $(u,v)\in U_M$ with $\mathrm{d}_M(u,v)< \tfrac{\rho}{14}$.
    \item  For $(u,w)\in M\times M$ with $\mathrm{d}_M(u,w)\geq \tfrac{\rho}{2}$ one has $h(u,v)(w) = w$.
    \item For all $u\in M$ we have $h(u,u) = \mathrm{id}_M$.
\end{itemize}

We define $\nu_D\colon E^f\to U_{D^f}$ as follows.
For $( \gamma,(x_0,x_1,\sigma),v)\in E^f $ set $$\nu_D(\gamma,(x_0,x_1,\sigma),v) = (\widetilde{\gamma},(x_0,x_1,\sigma))  \in U_{D^f}  $$ where
$$  \widetilde{\gamma}(t) =  h(\sigma(0),\exp_{\sigma(0)}(v))(\gamma(t)) \quad \text{for}\,\,\,t\in I .     $$
We define an explicit inverse by $\kappa_D\colon U_{D^f}\to E^f$ as follows.
Let $(\gamma,(x_0,x_1,\sigma))\in U_{D^f}$ and define 
$$\kappa_D (\gamma,(x_0,x_1,\sigma)) =  (\omega,(x_0,x_1,\sigma),v) \in E^f   $$
with
$$   \omega (t) =  h(\sigma(0),\gamma(0))^{-1}(\gamma(t)) \quad \text{and}\quad v = \exp_{\sigma(0)}^{-1}(\gamma(0)) .    $$
One checks that these maps are inverse to each other.
    The commutativity of the diagram in the second part of the lemma is then easy to see.
    Finally, one can use the ideas of the proof of \cite[Lemma 2.3]{hingston:2017} to show the last part.
\end{proof}

There is a map $d\colon D^f \to P^f$ given by
$$   d(\gamma,(p,q,\sigma)) = (p,q, \mathrm{concat}(\gamma,\sigma)) \quad \text{for}\,\,\, (\gamma,(p,q,\sigma))\in D^f .        $$

\begin{definition}
    Let $M$ be a closed oriented manifold and $N$ a closed manifold with a smooth map $f\colon N\to M$.
    We define the pairing $\cdot \colon \mathrm{H}_i(\Lambda M)\otimes \mathrm{H}_j(P^f) \to \mathrm{H}_{i+j-n}(P^f)$ as the composition
    \begin{eqnarray*}
	\cdot \colon  \mathrm{H}_{i}(\Lambda)\otimes \mathrm{H}_j(P^f) 
    &
    \xrightarrow{  (-1)^{n-ni}\times   }& \mathrm{H}_{i+j}(\Lambda\times P^f)
	\\ 
    &   \xrightarrow{ \hphantom{concii}\hphantom{conci} }& \mathrm{H}_{i+j}(\Lambda\times P^f, \Lambda\times P^f \setminus D^f) \\
	&\xrightarrow{\hphantom{ci} \text{excision}\hphantom{c}}  &
	\mathrm{H}_{i+j}(U_{D^f}, U_{D^f}\setminus D^f) 
    \\
	&\xrightarrow{\hphantom{coi} \tau_{D^f}\cap \hphantom{ic} }& \mathrm{H}_{i+j-n}(U_{D^f})
	\\
    &\xrightarrow[]{\hphantom{ci}{(\mathrm{R}_{D^f})_*} \hphantom{c}} & \mathrm{H}_{i+j-n}(D^f) 
	\\
    & \xrightarrow[]{\hphantom{iciii}d_*\hphantom{iiici}} & \mathrm{H}_{i+j-n}(P^f) \, .
\end{eqnarray*}

\end{definition}
\begin{remark}
    We remark that in the case that $f\colon N\hookrightarrow M$ is the inclusion of a submanifold, the space $D^f$ can be identified with the space
    $$  D_N =  \{ (\gamma,\sigma)\in \Lambda M\times P_N M\,|\, \gamma(0) = \sigma(0) \}  .   $$
    We shall use the notation $D_N$ in subsequent sections.
\end{remark}

\begin{prop}
    Let $M$ be a closed oriented manifold and let $N$ be a closed manifold.
    Assume that $f\colon N\to M$ is a smooth map and take homology with coefficients in a commutative ring $R$.
    Then the homology $\mathrm{H}_{\bullet}(P^f)$ is a left module over the Chas-Sullivan ring $(\mathrm{H}_{\bullet}(\Lambda M),\wedge_{\mathrm{CS}})$ via the pairing $\cdot \colon \mathrm{H}_{\bullet}(\Lambda M)\otimes \mathrm{H}_{\bullet}(P^f) \to \mathrm{H}_{\bullet-n}(P^f)$.
\end{prop}
The proof of this proposition uses very similar arguments to the proof of associativity and unitality of the Chas-Sullivan product, see e.g. \cite[Appendix B]{hingston:2017}.
We also used arguments along these lines in the proof of Theorem \ref{theorem_properties_path_product}.1.
Therefore we omit the proof of this proposition.

\section{The module structure in particular cases}\label{sec_module_structure_explicit}

Before we show the invariance of the module structure under homotopies of $f\colon N\to M$ in the next section we want to study some particular cases for the map $f\colon N\to M$.

\subsection{The constant map}

Assume that $f = \mathrm{cst}\colon N\to M$ is the constant map to the basepoint $p_0$ of $M$.
In this case we have
$$   P^{\mathrm{cst}} \cong N\times N\times \Omega M   \quad \text{and}\quad D^{\mathrm{cst}}  \cong \Omega M\times N\times N\times \Omega M . $$
Let $j\colon \Omega M\to \Lambda M$ be the inclusion of the based loop space and denote by $j_! \colon \mathrm{H}_{\bullet}(\Lambda M) \to \mathrm{H}_{\bullet- n}(\Omega M)$ the induced Gysin map.
We note that this is explicitly given as follows.
Let 
$$  U_{\Omega} =  \{ \gamma\in\Lambda M\,|\, \mathrm{d}_M(\gamma(0),p_0) <\epsilon'\} .   $$
Then $U_{\Omega}$ is a tubular neighborhood of $\Omega M$ in $\Lambda M$ and consequently, there is a Thom class $\tau_{\Omega}\in \mathrm{H}^n( U_{\Omega}, U_{\Omega}\setminus \Omega M)$.
The Gysin map $j_!$ is given by the composition
$$  \mathrm{H}_{\bullet}(\Lambda M) \xrightarrow[]{} \mathrm{H}_{\bullet}(U_{\Omega}, U_{\Omega }\setminus \Omega M)\xrightarrow[]{\tau_{\Omega}\cap} \mathrm{H}_{\bullet-n}(U_{\Omega}) \xrightarrow[]{(\mathrm{R}_{\Omega})_*} \mathrm{H}_{\bullet - n}(\Omega M) .      $$
Here, the first map is induced by the inclusion of pairs and excision and $\mathrm{R}_{\Omega}\colon U_{\Omega}\to \Omega M$ is the retraction map of the tubular neighborhood.

We now define a pairing $$\nu_{N,\Omega} \colon \mathrm{H}_i(\Lambda M ) \otimes \big(  \mathrm{H}_{\bullet}(N)\otimes \mathrm{H}_{\bullet}(N) \otimes \mathrm{H}_{\bullet}(\Omega M)\big)_j \to \big(  \mathrm{H}_{\bullet}(N)\otimes \mathrm{H}_{\bullet}(N) \otimes \mathrm{H}_{\bullet}(\Omega M)\big)_{i+j-n} . $$
For $A\in\mathrm{H}_{\bullet}(\Lambda M ),  a,b\in\mathrm{H}_{\bullet}(N),x\in\mathrm{H}_{\bullet}(\Omega M)$
we set
$$    \nu_{N,\Omega}( A\otimes a\otimes b\otimes x) =  (-1)^{(n+|A|)(n+|a|+|b|)} a\otimes b\otimes (j_! A \star x)   $$
where $\star$ is the Pontryagin product as before.
\begin{lemma}
    The pairing $\nu_{N,\Omega}$ makes $\mathrm{H}_{\bullet}(N)^{\otimes 2}\otimes \mathrm{H}_{\bullet}(\Omega)$ into a left module over the Chas-Sullivan ring.
\end{lemma}
\begin{proof}
    Let $A,B\in\mathrm{H}_{\bullet}(\Lambda M), a,b\in\mathrm{H}_{\bullet}(N)$ and $x\in\mathrm{H}_{\bullet}(\Omega M)$.
    We need to show that
    \begin{equation}\label{eq_module_structure_cst} \nu_{N,\Omega}( A\otimes \nu_{N,\Omega}(B\otimes a\otimes b\otimes x))  =  \nu_{N,\Omega} ( A\wedge_{\mathrm{CS}} B\otimes a\otimes b\otimes x) .      \end{equation}
    One sees directly that the left hand side is equal to the term
    $$     (-1)^{(|A|+|B|)(n+|a|+|b|)} a\otimes b\otimes (j_! A\star (j_!B \star x))  .        $$
    By \cite[Section 9.3]{goresky:2009}
    we have the identity
    $$     j_! A \star j_! B =  j_! (A\wedge_{\mathrm{CS}} B) \quad \text{for} \,\,\, A,B\in\mathrm{H}_{\bullet}(\Lambda M) .  $$
    For a detailed proof of this identity we refer to \cite[Section 2.3]{kupper2020homology}.
    Using this identity we obtain
    $$      \nu_{N,\Omega}( A\otimes \nu_{N,\Omega}(B\otimes a\otimes b\otimes x)) =     (-1)^{(|A|+|B|)(n+|a|+|b|)} a\otimes b\otimes (j_! (A\wedge_{\mathrm{CS}} B) \star x)  .     $$
    On the other hand, since $|A\wedge_{\mathrm{CS}} B| = |A | + |B| - n$, we see that the right hand side of equation \eqref{eq_module_structure_cst} yields the same expression.

    Recall that the unit of the Chas-Sullivan ring is given as follows.
    Let $s\colon M\to \Lambda M$ be the map 
    $$   s(p) = c_p\in \Lambda M    \quad \text{for}\,\,\,p\in M $$
    where $c_p\in \Lambda M$ is the constant loop at the point $p\in M$.
    If $[M]\in\mathrm{H}_n(M)$ is the fundamental class of $M$ then $s_*[M]\in\mathrm{H}_n(\Lambda M)$ is the unit of the Chas-Sullivan ring.
    Moreover we have
    $$   j_! (s_*[M]) =  [p_0] \in \mathrm{H}_0(\Omega M)     $$
    where $[p_0]$ is the generator of $\mathrm{H}_0(\Omega_0)$ with $\Omega_0\subseteq \Omega M$ being the component of the trivial loop.
    It is well-known that $[p_0]$ is the unit of the Pontryagin ring and consequently we see that
    $$   \nu_{N,\Omega}( s_*[M] \otimes a\otimes b\otimes x) =  a\otimes b\otimes x   $$
    for $a,b\in\mathrm{H}_{\bullet}(N)$ and $x\in\mathrm{H}_{\bullet}(\Omega M)$.
    This shows that 
    $$    \nu_{N,\Omega}  \colon  \mathrm{H}_{\bullet}(\Lambda)\otimes  \mathrm{H}_{\bullet}(N)^{\otimes 2}\otimes \mathrm{H}_{\bullet}(\Omega) \to \mathrm{H}_{\bullet}(N)^{\otimes 2}\otimes \mathrm{H}_{\bullet}(\Omega)   $$
    is a left module structure.
\end{proof}

\begin{prop}\label{prop_trivial_map_module}
    Let $M$ and $N$ be closed manifolds and assume that $M$ is oriented.
    Let $f\colon N\to M$ be the constant map and take homology with coefficients in a field.
    Then there is a commutative diagram
    $$
    \begin{tikzcd}
        \mathrm{H}_i(\Lambda M)\otimes \mathrm{H}_j(P^{f}) \arrow[]{r}{\cong} \arrow[]{d}{\cdot} 
        &
        [2.5em]
        \mathrm{H}_i(\Lambda M)\otimes \big( \mathrm{H}_{\bullet}(N)\otimes \mathrm{H}_{\bullet}(N)\otimes \mathrm{H}_{\bullet}(\Omega M)\big)_j \arrow[]{d}{\nu_{N,\Omega}}
        \\
        \mathrm{H}_{i+j-n}(P^f) \arrow[]{r}{\cong} &
        \big(  \mathrm{H}_{\bullet}(N)\otimes \mathrm{H}_{\bullet}(N)\otimes \mathrm{H}_{\bullet}(\Omega M)\big)_{i+j-n} .
    \end{tikzcd}
    $$
\end{prop}
\begin{proof}
    We shall write $\Lambda = \Lambda M$ and $\Omega = \Omega M$ throughout this proof.
    We will use the notation $(X,\sim A)$ for the topological pair $(X,X\setminus A)$.
    
    Recall from above that
    $$   P^{\mathrm{cst}} \cong N\times N\times \Omega    \quad \text{and}\quad D^{\mathrm{cst}}  \cong \Omega  \times N\times N\times \Omega  . $$
    Let $T\colon \Lambda \times P^{\mathrm{cst}}\to N\times N\times \Lambda \times \Omega$ be the obvious swapping map.
    We begin by claiming that the following diagram commutes.
    $$
    \begin{tikzcd}
        \mathrm{H}_i(\Lambda \times P^{\mathrm{cst}}) \arrow[]{r}{T_*}
        \arrow[]{d}{}
        & [3em]
        \mathrm{H}_i(N\times N\times \Lambda \times \Omega )\arrow[]{d}{}
        \\
        \mathrm{H}_i(U_{D^{\mathrm{cst}}},\sim  D^{\mathrm{cst}} ) \arrow[]{r}{T_*}
        \arrow[]{d}{\tau_{D^{\mathrm{cst}}}\cap}
        &
    \mathrm{H}_i(N\times N\times U_{\Omega} \times \Omega ,\sim N\times N\times \Omega \times \Omega )
    \arrow[]{d}{\mathrm{pr}_3^* \tau_{\Omega}\cap}
    \\
    \mathrm{H}_{i-n}(U_{D^{\mathrm{cst}}}) \arrow[]{r}{T_*}
    \arrow[]{d}{(d\circ \mathrm{R}_{D^{\mathrm{cst}}})_*}
    &
    \mathrm{H}_{i-n}(N\times N\times U_{\Omega}\times \Omega) \arrow[]{d}{h_*}
    \\
    \mathrm{H}_{i-n}(P^{\mathrm{cst}}) \arrow[]{r}{ = }
    &
    \mathrm{H}_{i-n}(N\times N\times \Omega ) .
    \end{tikzcd}
    $$
    Here, $h\colon N\times N\times U_{\Omega}\times \Omega \to N\times N\times \Omega$ is the map
    $$   h= \mathrm{id}_N\times \mathrm{id}_N\times (\mathrm{concat}\circ (\mathrm{R}_{\Omega}\times \mathrm{id}_{\Omega})) . $$
    Moreover, $\mathrm{pr}_3\colon N\times N\times U_{\Omega }\times \Omega \to U_{\Omega}$ is the projection on the third factor.
    The commutativity of the first and the last square is easy to see.
    For the second square, we have to show that
    $$   \mathrm{pr}_3^*\tau_{\Omega}\cap (T_* X)  =  T_* (\tau_{\mathrm{D}^{\mathrm{cst}}} \cap X)       $$
    for $X\in \mathrm{H}_i(U_{D^{\mathrm{cst}}}, \sim D^{\mathrm{cst}})$.
    Note that by naturality the left hand side of the above equation equals
    $$      T_* \big( ( T^* \mathrm{pr}_3^* \tau_{\Omega})\cap X\big)  . $$
    One checks directly that $T^* \mathrm{pr}_3^* \tau_{\Omega} = \tau_{D^{\mathrm{cst}}}$ and thus the second square commutes as well.
    Let $A\in\mathrm{H}_{\bullet}(\Lambda )$ and $a,b\in\mathrm{H}_{\bullet}(N)$ as well as $x\in\mathrm{H}_{\bullet}(\Omega M)$.
    Then the above shows that
    \begin{equation}\label{eq_module_signs_figured_out}
          A \cdot (a\times b\times x)  =  (-1)^{n-n|A|} \big( h_* (\mathrm{pr}_3^*\tau_{\Omega}\cap (T_* (A\times a\times b\times x) ) ) \big) .      
    \end{equation}   
    In order to complete the proof we consider the following diagram.
    $$
    \begin{tikzcd}
        \mathrm{H}_{\bullet}(\Lambda \times N^2\times \Omega) \arrow[]{d}{T_*}
        & [2.5em]
        \mathrm{H}_{\bullet}(\Lambda )\otimes \mathrm{H}_{\bullet}(N)^{\otimes 2} \otimes \mathrm{H}_{\bullet}(\Omega ) \arrow[swap]{l}{\times}
        \arrow[]{d}{\mathcal{T}}
        \\
        \mathrm{H}_{\bullet}(N^2\times \Lambda \times \Omega) \arrow[]{d}{\mathrm{pr}_3^*\tau_{\Omega}\cap}
        &
        \mathrm{H}_{\bullet}(N)^{\otimes 2} \otimes \mathrm{H}_{\bullet}(\Lambda) \otimes \mathrm{H}_{\bullet}(\Omega) \arrow[swap]{l}{\times}
        \arrow[]{d}{\mathrm{id}\otimes j_! \otimes \mathrm{id}}
        \\
        \mathrm{H}_{\bullet}(N^2\times \Omega \times \Omega) \arrow[]{d}{ (\mathrm{id}_{N^2}\times \mathrm{concat})_*}
        &
        \mathrm{H}_{\bullet}(N)^{\otimes 2}\otimes \mathrm{H}_{\bullet}(\Omega )^{\otimes 2}
        \arrow[swap]{l}{\times}
        \arrow[]{d}{\mathrm{id}\otimes \star} 
        \\
        \mathrm{H}_{\bullet}(N^2\times \Omega) & 
        \mathrm{H}_{\bullet}(N)^{\otimes 2}\otimes \mathrm{H}_{\bullet}(\Omega ) \arrow[swap]{l}{\times } .
    \end{tikzcd}
    $$
    Here, $$\mathcal{T}  \colon \mathrm{H}_{\bullet}(\Lambda )\otimes \mathrm{H}_{\bullet}(N)^{\otimes 2}\otimes \mathrm{H}_{\bullet}(\Omega ) \to  \mathrm{H}_{\bullet}(N)^{\otimes 2}\otimes  \mathrm{H}_{\bullet}(\Lambda ) \otimes \mathrm{H}_{\bullet}(\Omega )  $$
    is the map swapping the tensor factors.
    Let $A\in \mathrm{H}_{\bullet}(\Lambda), a,b,\in\mathrm{H}_{\bullet}(N)$ and $x\in\mathrm{H}_{\bullet}(\Omega )$.
    Then the first square commutes up to sign $(-1)^{(|a|+|b|)|A|}$ and the second one up to sign $(-1)^{n(|a|+|b|)}$.
    Hence, we get that
    $$  h_* (\mathrm{pr}_3^*\tau_{\Omega}\cap (T_*( A\times a\times b\times x ) ))  =   (-1)^{(n+|A|)(|a|+|b|)} a\otimes b\otimes (j_! A\star x) .   $$
    Combining this with equation \eqref{eq_module_signs_figured_out} yields the claim.
\end{proof}

Once we have shown that the module structure on $\mathrm{H}_{\bullet}(P^f)$ over the Chas-Sullivan ring is invariant under homotopies of $f$, we will see that the above result yields the module structure in the case that the map $f\colon N\to M$ is null-homotopic.

\subsection{Embedding of a submanifold}

We now turn to the situation that $f\colon N\to M$ is the embedding of a submanifold.
In this situation the module pairing factors in a particular way as we shall now explain.

Let $M$ be a closed oriented manifold and $N$ a closed submanifold of $M$.
Let $\Lambda_N M\subseteq \Lambda M$ be the subspace
$$  \Lambda_N M = \{ \gamma\in \Lambda M\,|\, \gamma(0) \in N\} .     $$
This is a codimension $r$ submanifold of $\Lambda M$ with $r = \mathrm{codim}(N\hookrightarrow M)$.
We denote the inclusion by $j\colon \Lambda_N M\hookrightarrow \Lambda M$.
This induces a Gysin map 
$$   j_!\colon \mathrm{H}_i(\Lambda M) \to \mathrm{H}_{i-r}(\Lambda_N M) .   $$
We now define a pairing
$$    \triangleright \colon \mathrm{H}_i(\Lambda_N M) \otimes \mathrm{H}_l(P_N M)\to \mathrm{H}_{i+l-r}(P_N M) .   $$
Let $D_N' = D_N \cap (\Lambda_N M \times P_N M)$, i.e. we have
$$  D_N' = \{ (\gamma,\sigma)\in \Lambda_N M\times P_N M\,|\, \gamma(0) = \sigma(0)\} .        $$
Note that $D_N\cong D_N'$ but we use the prime to make clear that we consider $D_N'$ as a subspace of $\Lambda_N M\times P_N M$.
Note that $D_N'$ is the preimage 
$$   D_N' = \big((\mathrm{ev}_{\Lambda}\times \mathrm{ev}_0 )|_{\Lambda_N M\times P_N M} \big)^{-1} (\Delta N) .   $$
Let $U_N$ be the tubular neighborhood of $\Delta N$ in $N\times N$ as before and set
$$   U_{D_N'} = \big((\mathrm{ev}_{\Lambda}\times \mathrm{ev}_0 )|_{\Lambda_N M\times P_N M} \big)^{-1}(U_N) .    $$
Then the Thom class $\tau_N\in\mathrm{H}^k(U_N,U_N\setminus N)$ pulls back to a Thom class
$$     \tau_{D_N'} = (\mathrm{ev}_{\Lambda}\times \mathrm{ev}_0)^* \tau_{N} \in \mathrm{H}^k(U_{D_N'},U_{D_N'}\setminus D_N')   .  $$
Moreover, there is a retraction map $\mathrm{R}_{D_N'}\colon U_{D_N'}\to D_N'$ and the concatenation map restricts to a map $\mathrm{concat}\colon D_N'\to P_N M$.
We note that the map $R_{D_N'}$ can be chosen as follows.
By definition we have
$$    U_{D_N'} = \{(\gamma,\sigma)\in \Lambda_N M\times P_N M\,|\, \mathrm{d}(\gamma(0),\sigma(0)) < \epsilon' \}   .   $$
We define
$$   \mathrm{R}_{D_N'} (\gamma,\sigma) = \big(  \mathrm{concat}(\widehat{\sigma(0)\gamma(0)},\gamma, \widehat{\gamma(0)\sigma(0)}),\,\sigma\big) \quad \text{for}\,\,\,(\gamma,\sigma)\in U_{D_N'} .    $$
As before the notation $\widehat{pq}$ means the unique length-minimizing geodesic between points $p,q\in M$ with $\mathrm{d}_M(p,q)<\epsilon'$.
With the same ideas as in the proof of Lemma \ref{lemma_retraction_explicit} one can see that $\mathrm{R}_{D_N'}$ is homotopic to the retraction induced by the normal bundle of $D_N'$.

\begin{definition}\label{def_triangleright}
	Let $M$ be a closed oriented manifold and let $N$ be a closed oriented submanifold of $M$.
	Take homology with coefficients in a commutative ring $R$.
	Then we define a pairing $\triangleright\colon\mathrm{H}_i(\Lambda_N M) \otimes \mathrm{H}_j(P_N M)\to \mathrm{H}_{i+j-k}(P_N M)  $ as the composition 
	\begin{eqnarray*}
		\triangleright\colon \mathrm{H}_i(\Lambda_N M)  \otimes \mathrm{H}_j(P_N M) & \xrightarrow[]{(-1)^{k-ki}\times} &
		\mathrm{H}_{i+j}(\Lambda_N M\times P_N M) \\
		&\xrightarrow[]{\hphantom{coiiincoccc}} & \mathrm{H}_{i+j}(\Lambda_N M\times P_N M,\Lambda M\times P_N M \setminus D_N') \\
		&\xrightarrow[]{\hphantom{ii}\mathrm{excision}\hphantom{ii}} &    \mathrm{H}_{i+j}(U_{D_N'},U_{D_N'}\setminus  {D_N'}) \\
		& \xrightarrow[]{\hphantom{cii}\tau_{D_N'}\cap\hphantom{cii}} & \mathrm{H}_{i+j-k}(U_{D_N'}) \\
		& \xrightarrow[]{\hphantom{ci}(R_{D_N'})_*\hphantom{c}}  & \mathrm{H}_{i+j-k}(D_N') \\
		& \xrightarrow[]{\hphantom{i}(\mathrm{concat})_*\hphantom{i}} & \mathrm{H}_{i+j-k}(P_N M) .    
	\end{eqnarray*}
\end{definition}

The next proposition relates the pairing $\triangleright$ to the module structure of $\mathrm{H}_{\bullet}(P_N M)$ over the Chas-Sullivan ring.

\begin{prop}\label{module_factors}
	Let $M$ be a closed oriented manifold and let $N$ be a closed oriented submanifold.
	Let $k$ be the dimension of $N$ and $r$ be the codimension.
	Let $j\colon \Lambda_N M\hookrightarrow \Lambda M$ be the inclusion and $\triangleright\colon\mathrm{H}_i(\Lambda_N M) \otimes \mathrm{H}_l(P_N M)\to \mathrm{H}_{i+l-r}(P_N M) $ be the pairing as in the above definition.
	Then the module structure $\cdot \colon \mathrm{H}_i(\Lambda M)\otimes \mathrm{H}_l(P_N M)\to \mathrm{H}_{i+l-n}(P_N M)$ is equal to the composition
	$$           \mathrm{H}_i(\Lambda M)\otimes \mathrm{H}_l(P_N M) \xrightarrow[]{(-1)^{r-ri+kr} \, j_! \otimes \id} \mathrm{H}_{i-r}(\Lambda_N M)\otimes \mathrm{H}_l(P_N M) \xrightarrow[]{\hphantom{co}\triangleright\hphantom{co}} \mathrm{H}_{i+l-n}(P_N M) .         $$   
\end{prop}
\begin{proof}
	The Gysin map $j_!\colon \mathrm{H}_i(\Lambda M)\to \mathrm{H}_{i-r}(\Lambda_N M)$ can be defined as the composition
	$$     \mathrm{H}_i(\Lambda M) \to \mathrm{H}_i(\Lambda M,\Lambda M\setminus \Lambda_N M) \to \mathrm{H}_i(\overline{U},\overline{U}\setminus \Lambda_N M) \xrightarrow[]{\overline{\tau}\cap} \mathrm{H}_{i-r}(\overline{U}) \xrightarrow[]{\rho_*} \mathrm{H}_{i-r}(\Lambda_N M) .     $$
	Here, $\overline{U}$ is a tubular neighborhood of $\Lambda_N M$ in $\Lambda M$, $\overline{\tau}\in \mathrm{H}^r(\overline{U},\overline{U}\setminus \Lambda_N M)$ is the corresponding Thom class and $\rho\colon \overline{U}\to \Lambda_N M$ is the retraction.
	In particular, $\overline{U}$ is homeomorphic to the normal bundle $E$ of $\Lambda_N M$ in $\Lambda M$ which is isomorphic to the bundle $E = \mathrm{ev}_{\Lambda}^* \mathcal{N}(N\hookrightarrow M)$.
	Recall that 
	$$     \mathcal{N}(\Delta N \hookrightarrow M\times N ) \cong TN\oplus \mathcal{N}(N\hookrightarrow M)   \cong i^* TM    $$
	where $\mathcal{N}(N\hookrightarrow M)$ is the normal bundle of the inclusion $i\colon N\hookrightarrow M$.
	Moreover, note that both $TN$ and $TM$ already carry an orientation so the orientation of $\mathcal{N}(N\hookrightarrow M)$ is chosen so that the above isomorphism preserves orientations.
	Therefore the normal bundle of $D\subseteq \Lambda M\times P_N M$ is isomorphic to
	$$          (\mathrm{ev}_{\Lambda}\times \mathrm{ev}_0)^*  \mathcal{N}(\Delta N \hookrightarrow M\times N ) \cong   (\mathrm{ev}_{\Lambda}\times \mathrm{ev}_0)^* TN \oplus (\mathrm{ev}_{\Lambda}\times \mathrm{ev}_0)^* (\mathcal{N}(N\hookrightarrow M))    $$
	where we understand $(\mathrm{ev}_{\Lambda}\times \mathrm{ev}_0) \colon D_N'\to \Delta N$.
	Let $\pi_1\colon E_1\to D_N'$ be the bundle 
    $$(\mathrm{ev}_{\Lambda}\times \mathrm{ev}_0)^* TN\to D_N' $$ 
    and $\pi_2\colon E_2\to D_N'$ be the bundle $$(\mathrm{ev}_{\Lambda}\times \mathrm{ev}_0)^* (\mathcal{N}(N\hookrightarrow M)) \to D_N'.$$
	Then it is well-known that the Thom class of $E_1\oplus E_2\to D_N'$ is equal to the cup-product
	$$      \pi_2^*(\tau_1) \cup \pi_1^*(\tau_2) \in \mathrm{H}^n(E,E\setminus D_N' )       $$ 
	where $\tau_1\in \mathrm{H}^k(E_1,E_1\setminus D_N')$ and $\tau_2\in \mathrm{H}^r(E_2,E_2\setminus D_N')$.
	By the identity relating the cap and the cup product, see \cite[Theorem VI.5.2]{bredon:2013} we see that for $X\in\mathrm{H}_i(\Lambda M\times P_N M)$ it holds that
    \begin{equation} \label{eq_gysin_commutes}
    (\mathrm{R}_{D_N})_*(\tau_{D_N}\cap X ) =   (\mathrm{R}_{D_N'})_* \big(\tau_{D_N'}\cap (\rho\times \mathrm{id})_*(\mathrm{pr}_1^*\overline{\tau}\cap X)\big) .  \end{equation}
	Now, let $X\in\mathrm{H}_{\bullet}(\Lambda M)$ and $Y\in\mathrm{H}_{\bullet}(P_N M)$.
	Then by equation \eqref{eq_gysin_commutes} we see that 
	$$     X\cdot Y \quad \text{and} \quad( j_! X )\triangleright Y   $$
	agree up to the sign $(-1)^{r-r|X| - kr}$.
	This completes the proof.
\end{proof}
\begin{remark}\label{remark_lie_groups_gysin}
	As we shall see in the following the above proposition is useful since in certain cases the Gysin map $j_!\colon \mathrm{H}_i(\Lambda M)\to \mathrm{H}_{i-r}(\Lambda_N M)$ can be understood nicely.
	If $G$ is a compact Lie group then there is a homeomorphism $\widetilde{F}\colon \Lambda G\to G\times \Omega_e G$ given by
	$$   \widetilde{F}(\gamma) = ( \gamma(0),\gamma(0)^{-1}\cdot \gamma) .      $$
	If $N$ is a closed submanifold of $G$ then one checks that $\widetilde{F}$ restricts to a homeomorphism $F\colon \Lambda_N G\to N\times \Omega_e G$.
	We obtain a commutative diagram
	$$     \begin{tikzcd}
		\Lambda_N G\arrow[]{r}{j} \arrow[swap]{d}{F} & [2.em] \Lambda G \arrow[]{d}{\widetilde{F}} 
		\\
		N\times \Omega G \arrow[]{r}{i\times \mathrm{id}_{\Omega}} & G\times \Omega G 
	\end{tikzcd}      $$
	where $i\colon N\hookrightarrow G$ is the inclusion.
	Consequently, the Gysin map $j_!$ can be computed using the Gysin map $i_!\colon \mathrm{H}_l(G)\to \mathrm{H}_{l-r}(N)$.
	Since the map $i$ is a map between finite-dimensional manifolds, we have much better computational tools available for the Gysin map $i_!$ compared to the Gysin map $j_!$.
\end{remark}

\subsection{The identity map}

We end this section by considering the identity map $\mathrm{id}\colon M\to M$ for a closed oriented manifold $M$.
Recall from Section \ref{sec_path_product_examples} that the map $\varphi_0\colon P^{\mathrm{id}} = PM \to M$ given by $\varphi_0(\gamma) = \gamma(0)$ is a homotopy equivalence.
\begin{prop}
    Let $M$ be a closed oriented manifold of dimension $n$ and consider the identity map $\mathrm{id}\colon M\to M$.
    Then the module structure $\mathrm{H}_{\bullet}(\Lambda M)\otimes \mathrm{H}_{\bullet}(PM) \to \mathrm{H}_{\bullet}(PM)$ is given by the augmentation of the Chas-Sullivan ring over the intersection ring of $M$. More precisely, the following diagram commutes
    $$    
    \begin{tikzcd}
        \mathrm{H}_i(\Lambda M)\otimes \mathrm{H}_j(PM) \arrow[]{r}{\cdot} 
        \arrow[d, "\mathrm{id}\otimes (\varphi_0)_*", "\cong"']
        & [3em]
        \mathrm{H}_{i+j-n}(P M) \arrow[d, "(\varphi_0)_*" , "\cong"'] 
        \\
        \mathrm{H}_i(\Lambda M)\otimes \mathrm{H}_j(M) \arrow[]{r}{(\mathrm{ev}_{\Lambda})_* ( \cdot) \overline{\cap} (\cdot)} 
        &
        \mathrm{H}_{i+j-n}(M) .
    \end{tikzcd}
    $$
\end{prop}
\begin{proof}
    We begin by showing that the following diagram commutes.
     $$
        \begin{tikzcd}
            \mathrm{H}_{i}(\Lambda M) \otimes \mathrm{H}_{j}(P^{\mathrm{id}}) \arrow[]{r}{ (\mathrm{ev}_{\Lambda})_*\otimes (\varphi_0)*}
            \arrow[]{d}{}
            & [2.5em]
            \mathrm{H}_{i}( M)\otimes \mathrm{H}_{j}(M) \arrow[]{d}{} 
            \\
            \mathrm{H}_{i+j}(U_{C^{\mathrm{id}}}, U_{C^{\mathrm{id}}} \setminus C^{\mathrm{id}}) \arrow[]{r}{(\mathrm{ev}_{\Lambda}\times \varphi_0)_*}
            \arrow[]{d}{\tau_{C^{\mathrm{id}}}\cap}
            & 
            \mathrm{H}_{i+j}(U_{M} , U_M  \setminus \Delta M)\arrow[]{d}{\tau_{M}\cap}
            \\
            \mathrm{H}_{i+j-n}(U_{C^{\mathrm{id}}}) \arrow[]{r}{(\mathrm{ev}_{\Lambda}\times \varphi_0)_*}
            \arrow[]{d}{(\mathrm{R}_{C^{\mathrm{id}}})_*}
            &
            \mathrm{H}_{i+j-n}(U_M)\arrow[]{d}{( \mathrm{R}_{M})_*}
            \\
            \mathrm{H}_{i+j-n}(C^{\mathrm{id}}) \arrow[]{r}{(\mathrm{ev}_{\Lambda}\times \varphi_0)_*}  \arrow[]{d}{c_*} &
            \mathrm{H}_{i+j-n}(M) 
            \\
            \mathrm{H}_{i+j-n}(P^{\mathrm{id}}) \arrow[swap]{ur}{(\varphi_0)_*} & 
        \end{tikzcd}
        $$
        The commutativity of the first square is easy to see.
        The second square commutes by naturality.
        Indeed, let $X\in\mathrm{H}_{\bullet}(U_C^{\mathrm{id}}, U_{C^{\mathrm{id}}}\setminus C^{\mathrm{id}})$, then we have
        $$    \tau_M\cap (\mathrm{ev}_{\Lambda}\times \varphi_0 )_* X = (\mathrm{ev}_{\Lambda}\times \varphi_0 )_* \big( (\mathrm{ev}_{\Lambda}\times \varphi_0 )^* \tau_{M} \cap X \big)           $$
        by naturality of the cap product.
        By definition we have $\varphi_0 = \mathrm{ev}_0$ and we recall that $\tau_{C^{\mathrm{id}} } = (\mathrm{ev}_{\Lambda}\times \varphi_0)^*\tau_M$.
        This shows that the second square commutes.
        The commutativity of the third square and of the lower triangle follows directly from strict commutativity of the underlying diagram of maps.
        Recall that the signs for the module structure and for the intersection product on $M$ are the same and consequently we obtain for $A\in \mathrm{H}_{i}(\Lambda M)$ and $X\in\mathrm{H}_j(PM)$ that
        $$   (\mathrm{ev}_{\Lambda})_* (A) \, \overline{\cap} \, (\varphi_0)_* (X) =   (\varphi_0)_* (A\cdot X)     . $$
        This shows the claim of the proposition.
\end{proof}

\section{Invariance of the module structure}\label{sec_invariance_module}

We now prove that the module structure on $\mathrm{H}_{\bullet}(P^f)$ is invariant under homotopies of $f$.
We will start by constructing a homotopy equivalence $\Lambda M\times P^f\simeq \Lambda M\times P^g$ which preserves $D^f$ and its tubular neighborhood.

Let $M$ and $N$ be closed oriented manifolds and let $f,g\colon N\to M$ be smooth maps.
In this section we shall assume without loss of generality that the tubular neighborhood of the diagonal $\Delta M\subseteq M\times M$ which is given by
$$  U_M =  \{(p,q)\in M\times M\,|\, \mathrm{d}_M(p,q) <\epsilon' \}  $$
is such that $3\epsilon'$ is smaller than the injectivity radius of $M$.
Recall that the pullback bundle $f^* TM \to N$ is given by
$$  f^* TM  = \{ (x,v)\in N\times TM\,|\, v\in T_{f(x)} M\} .     $$
Moreover, recall that we have an explicit tubular neighborhood of the diagonal $\Delta M \subseteq U_M$ given by
$$  \varphi\colon TM\to U_M ,\quad (p,v) \mapsto (p,\exp_p(v))        $$
where we identify $TM$ with the open disk bundle $\{v\in TM\,|\, |v|<\epsilon'\}$.
Define a map $\pi_f\colon U_{D^f}\to f^* TM$ by
$$    \pi_f (\gamma,(x_0,x_1,\sigma))  =  (x_0,(\exp_{f(x_0)})^{-1} (\gamma(0))) .    $$
By definition of $U_{D^f}$ we have $\mathrm{d}_M( f(x_0),\gamma(0))< \epsilon'$, hence this is well-defined.
Note that the map $\mathrm{ev}_{\Lambda }\times \mathrm{ev}_0\colon U_{D^f}\to U_M$ factors through $\pi_f$ which can be expressed by the commutativity of the following diagram.
$$ 
\begin{tikzcd}
    U_{D^f}\arrow[]{r}{\mathrm{ev}_{\Lambda }\times \mathrm{ev}_0}
    \arrow[]{d}{\pi_f}
    & [2em]
    U_M
    \\
    f^* TM \arrow[]{r}{\mathrm{pr}_{TM}}
    &
    TM \arrow[swap]{u}{\varphi} .
\end{tikzcd}
$$
Here, $\mathrm{pr}_{TM}\colon f^* TM\to TM$ is the canonical map induced by the construction of $f^* TM$ as a pullback.
We also remark that $D^f =  \pi_f^{-1}(0)$ is the preimage of the zero-section $N\to f^* TM$.

\begin{lemma}\label{lemma_construction_homotopy_equivalence}
    Let $M$ and $N$ be closed oriented manifolds and let $f,g\colon N\to M$ be smooth homotopic maps.
    Fix a homotopy $H\colon N\times I\to M$ between $f$ and $g$ as well as a bundle isomorphism $F \colon f^* TM \to g^* TM$.
    Then there is a map $\Gamma\colon \Lambda M\times P^f \to \Lambda M\times P^g$ which satisfies the following properties:
    \begin{enumerate}
        \item The map $\Gamma$ is homotopic to $\mathrm{id}_{\Lambda M}\times \Phi\colon \Lambda M\times P^f\to \Lambda M\times P^g$ where $\Phi$ is the map described in the beginning of Section \ref{sec_invariance_product}.
        \item The diagram
        $$  
        \begin{tikzcd}
            U_{D^f} \arrow[]{r}{\Gamma}
            \arrow[]{d}{\pi_f} 
            &
            U_{D^g} \arrow[]{d}{\pi_g}
            \\
            f^* TM \arrow[]{r}{F} & g^* TM 
        \end{tikzcd}
        $$
        commutes.
        \item In particular, $\Gamma$ induces a map of pairs
        $   (U_{D^f}, U_{D^f}\setminus D^f  ) \to (U_{D^g}, U_{D^g}\setminus D^g)  .    $
    \end{enumerate}
\end{lemma}
    Note that the last property follows directly from the second property, we just spell it out for convenience.
\begin{proof}
    We shall construct the map $\Gamma\colon \Lambda M\times P^f\to \Lambda M\times P^g$ explicitly.
    Let $H\colon N\times I\to M$ be a fixed homotopy between $f$ and $g$ and let $\eta\colon N\to PM$ be the adjoint, i.e.
    $$  \eta(p)(s) =  H(p,s)   $$
    for all $p\in N$.
    In particular we have
    $$  \eta(p)(0) =  f(p) \quad \text{and}\quad \eta(p)(1) = g(p)  .  $$
    Define the set
    $$  \widetilde{U}_f = \{ (x,p)\in N\times M\,|\,  \mathrm{d}_M(f(x),p)< 3 \epsilon'  \}  .  $$
    Moreover, define $\xi\colon \widetilde{U}_f \to M$ by
    $$   \xi_{(x,p)} :=  \xi(x,p) : =   \exp_{g(x)}\big( 
    \mathrm{pr}_{TM}( F(x, \exp_{f(x)}^{-1}( p )) \big)    $$
    for $(x,p)\in\widetilde{U}_f $.
    As above, by $\mathrm{pr}_{TM}\colon g^* TM\to TM$ we mean the canonical map $g^* TM \to TM$ which is given by the pullback diagram.
    Note that the map $\xi$ is well-defined, since 
    $$     \exp_{f(x)}^{-1}(p) \in T_{f(x_0)} M    $$
    and thus 
    $$  F(x, \exp_{f(x)}^{-1}(p)) \in (g^* TM) \big|_{x}   .   $$
    By construction we have
    $    \mathrm{d}_M( g(x), \xi(x,p) ) < 3 \epsilon'     $.
    Define a map 
    $\Xi\colon   \widetilde{U}_f \to PM          $
    by setting
    $$    \Xi(x,p) =   \mathrm{concat}(  \widehat{\xi_{(x,p)} g(x)}, \overline{\eta(x)}, \widehat{f(x)p} ) .       $$
    Here $\widehat{qr}\in PM$ denotes the unique length-minimizing geodesic segment joining $q\in M$ and $r\in M$
    and $\overline{\omega}\in PM$ is the reverse path of the path $\omega\in PM$.
    We refer to Figure \ref{pic_homotopy} for a sketch of the path $\Xi(x,p)$.
    \begin{figure}[t]
\centering
\includegraphics[scale=0.4]{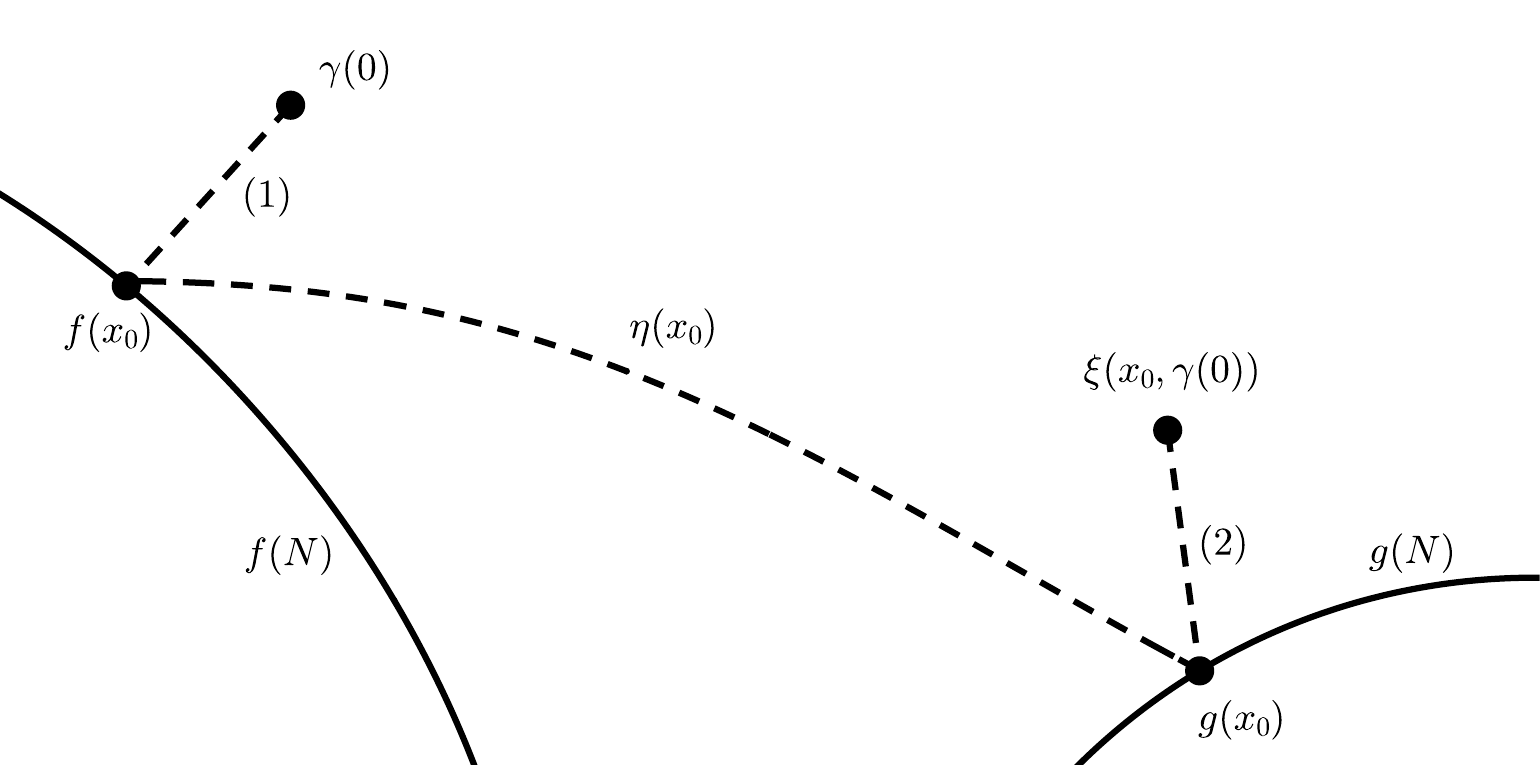}
\caption{The path $\Xi(x_0,\gamma(0))$. The dashed paths labelled (1) and (2) are the geodesic sticks from $\gamma(0)$ to $f(x_0)$ and from $g(x_0)$ to $\xi(x_0,\gamma(0))$, respectively. The dashed path labelled $\eta(x_0)$ is the path defined by the homotopy $H$ between $f(x_0)$ and $g(x_0)$.}
\label{pic_homotopy}
\end{figure}
    Let 
    $$   U_f =  \{ (\gamma,(x_0,x_1,\sigma))\in \Lambda M\times P^f\,|\, \mathrm{d}_M(\gamma(0),\sigma(0)) < 3\epsilon' \}      $$
    and analogously define $U_g$.
    We define a map 
    $   \widetilde{\Gamma} \colon U_{f}\to U_{g}      $
    by setting
    $$    \widetilde{\Gamma}((\gamma,(x_0,x_1,\sigma))) =   ( \gamma', \Phi(x_0,x_1,\sigma))      $$
    with $\Phi$ as in the beginning of Section \ref{sec_invariance_product} and with 
    \begin{equation}\label{eq_gamma_prime}
      \gamma' =  \mathrm{concat}(  \Xi(x_0,\gamma(0)) \,,\, \gamma \, ,\, \overline{\Xi(x_0,\gamma(0))}) .     
    \end{equation}
    In order to define $\Gamma\colon \Lambda M\times P^f\to \Lambda M\times P^g$ set $\Gamma = \widetilde{\Gamma}$ on $U_{D^f}$.
    On the set
    $$   W =    \{ (\gamma,(x_0,x_1,\sigma))\in \Lambda M\times P^f\,|\, \epsilon' \leq  \mathrm{d}_M(\gamma(0),\sigma(0)) \leq 2\epsilon' \} \subseteq \widetilde{U}_f    $$
    we define $\Gamma$ as follows.
    Let $\delta\colon \Lambda M\times P^f \to \mathbb{R}$ be the map $\delta(\gamma,(x_0,x_1,\sigma)) = \mathrm{d}_M(\gamma(0),\sigma(0))  $.
    Then define $\Gamma\colon W\to \Lambda \times P^g$ as follows.
    Recall that we have $\widetilde{\Gamma}(\gamma,(x_0,x_1,\sigma)) =  (\gamma',\Phi(x_0,x_1,\sigma))$.
    Then we set
    $$  \Gamma((\gamma,(x_0,x_1,\sigma))) =    ( \gamma' |_{[a,b]}, \Phi(x_0,x_1,\sigma))      $$
    where
    $$   a =  \frac{\delta(\gamma,(x_0,x_1,\sigma))}{3\epsilon} - \frac{1}{3} \quad \text{and}\quad b=   \frac{4}{3} - \frac{\delta(\gamma,(x_0,x_1,\sigma))}{3\epsilon}  .      $$
    In particular we have
    $$ a =  0 \quad \text{and} \quad b= 1 \quad \text{for} \quad \delta(\gamma,(x_0,x_1,\sigma)) = \epsilon   $$
    and
    $$   a = \frac{1}{3} \quad \text{and}\quad b= \frac{2}{3} \quad \text{for}    \quad \delta(\gamma,(x_0,x_1,\sigma)) =  2 \epsilon .  $$
    Note that all the restrictions of $\gamma'$ are indeed loops by the construction of $\gamma'$.
    Finally, define $\Gamma$ to be equal to $\mathrm{id}_{\Lambda M}\times \Phi$ on $\Lambda M\times P^f \setminus (U_{D^f}\cup W)$.
    The resulting map $\Gamma\colon \Lambda M\times P^f \to \Lambda M\times P^g$ is continuous.
    One checks from the construction that the claimed properties are satisfied.
\end{proof}

\begin{theorem}\label{theorem_invariance_module_structure}
    Let $M$ and $N$ be closed manifolds and assume that $M$ is oriented.
    Let $f,g\colon N\to M$ be smooth homotopic maps.
    Take homology with coefficients in a commutative ring $R$.
    Then the homotopy equivalence $\Phi\colon P^f\to P^g$ induces an isomorphism of the module structures of $\mathrm{H}_{\bullet}(P^f)$ and $\mathrm{H}_{\bullet}(P^g)$ over the Chas-Sullivan ring of $M$.
    More precisely, the diagram
    $$   
    \begin{tikzcd}
        \mathrm{H}_{i}(\Lambda M) \otimes \mathrm{H}_j(P^f) \arrow[]{r}{\mathrm{id}\otimes \Phi_*}
        \arrow[]{d}{\cdot}
        &
        [3em]
        \mathrm{H}_i(\Lambda M)\otimes \mathrm{H}_j(P^g) \arrow[]{d}{\cdot}
        \\
        \mathrm{H}_{i+j-n}(P^f) \arrow[]{r}{\Phi_*}
        &
        \mathrm{H}_{i+j-n}(P^g) 
    \end{tikzcd}
    $$
    commutes.
\end{theorem}
\begin{proof}
    As before we shall use the notation $\Lambda$ for $\Lambda M$ throughout this proof as well as the notation $(X,\sim A)$ for the topological pair $(X,X\setminus A)$.
    We claim that the following diagram commutes.
    $$
    \begin{tikzcd}
        \mathrm{H}_i(\Lambda \times P^f) \arrow[]{r}{(\mathrm{id}_{\Lambda}\times \Phi)_*}
        \arrow[]{d}{}
        & [2.5em]
        \mathrm{H}_i(\Lambda \times P^g) \arrow[]{d}{}
        \\
        \mathrm{H}_i(U_{D^f},\sim D^f) \arrow[]{r}{\Gamma_*} \arrow[swap]{d}{\tau_{D^f}\cap}
        &
        \mathrm{H}_i(U_{D^g},\sim D^g) \arrow[]{d}{\tau_{D^g}\cap}
        \\
        \mathrm{H}_{i-n}(U_{D^f}) \arrow[]{r}{\Gamma_*}
        \arrow[swap]{d}{(d\circ \mathrm{R}_{D^f})_*}&
        \mathrm{H}_{i-n}(U_{D^g}) \arrow[]{d}{(d\circ \mathrm{R}_{D^g})_*}
        \\
        \mathrm{H}_{i-n}(P^f) \arrow[]{r}{\Phi_*}
        &
        \mathrm{H}_{i-n}(P^g) .
    \end{tikzcd}
    $$
    The commutativity of the first square follows from Lemma \ref{lemma_construction_homotopy_equivalence}.
    For the second square let $X\in\mathrm{H}_i(U_{D^f},\sim D^f)$.
    We want to show that
    $$   \tau_{D^g}\cap (\Gamma_* X )   =   \Gamma_* (\tau_{D^f}\cap X)  .      $$
    By naturality the left hand side of this equation is equal to 
    $$   \Gamma_* (  (\Gamma^*\tau_{D^g})\cap X       )   .  $$
    Recall that $\tau_{D^g} = (\mathrm{ev}_{\Lambda}\times \mathrm{ev}_0)^* \tau_M$ and that we have
    $$      \mathrm{ev}_{\Lambda} \times \mathrm{ev}_0|_{U_{D^g}} =  \varphi\circ \mathrm{pr}_{TM}\circ \pi_g  .  $$
    In particular we have that
    $$\tau_g\colon = (\varphi\circ \mathrm{pr}_{TM})^* \tau_M  \in \mathrm{H}^n( g^* TM,\sim N)   $$ 
    is the Thom class of the bundle $g^* TM\to N$.
    Moreover by Lemma \ref{lemma_construction_homotopy_equivalence} we have
    $$    \pi_g\circ \Gamma =  F\circ \pi_f,     $$
    thus we get
    $$  \Gamma^* \tau_{D^g} =  \pi_f^* F^* \tau_g .   $$
    Clearly the bundle isomorphism $F\colon f^* TM\to g^* TM$ pulls the Thom class $\tau_g$ back to the Thom class
    $$    \tau_f\in \mathrm{H}^n(f^* TM,\sim N) .     $$
    Hence, we have
    $$   \Gamma^* \tau_{D^g} =  \pi_f^* \tau_f =  \tau_{D^f}     $$
    and this shows the commutativity of the second square.
    For the commutativity of the last square we claim that the underlying diagram of maps commutes up to homotopy.
    Indeed we compute that for $(\gamma,(x_0,x_1,\sigma))\in \Lambda M\times P^f$ we have that up to global re-parametrization
     $$     \Phi\circ d\circ \mathrm{R}_{D^f} (\gamma,(x_0,x_1,\sigma)) =  (x_0,x_1,\mathrm{concat}(\overline{\eta(x_0)},\widehat{\sigma_0\gamma_0},\gamma,\widehat{\gamma_0\sigma_0},\sigma,\eta(x_1)) ,     $$
      where we wrote $\gamma_0 = \gamma(0)$ and $\sigma_0 = \sigma(0)$.
     On the other hand we obtain that up to global re-parametrization
    \begin{align}
        \nonumber
        d\circ \mathrm{R}_{D^g}\circ \Gamma (\gamma,(x_0,x_1,\sigma)) &=& \big(x_0,x_1, \mathrm{concat}( 
        {\widehat{g(x_0)\xi_{(x_0,\gamma_0)}}, \widehat{\xi_{(x_0,\gamma_0)}g(x_0)}},  
        \overline{\eta(x_0)}, \widehat{\sigma_0\gamma_0}, \gamma,   \widehat{\gamma_0\sigma_0},
        \\ 
        & \hphantom{bla} &   \nonumber
        {\eta(x_0), \widehat{g(x_0)\xi_{(x_0,\gamma_0)}}\widehat{\xi_{(x_0,\gamma_0)}g(x_0)} , \overline{\eta(x_0)}},\sigma,\eta(x_1)) \big) 
    \end{align}
    Note that the concatenation of the paths 
    $$  \widehat{g(x_0)\xi_{(x_0,\gamma_0)}} \quad \text{and} \quad \widehat{\xi_{(x_0,\gamma_0)}g(x_0)}   $$
    as well as of the paths 
    $$  \eta(x_0), \widehat{g(x_0)\xi_{(x_0,\gamma_0)}}, \widehat{\xi_{(x_0,\gamma_0)}g(x_0)}   \quad \text{and}\quad \overline{\eta(x_0)}$$
    can be contracted by the "spaghetti trick".
    Indeed, we see that this yields a homotopy between $ d\circ \mathrm{R}_{D^g}\circ \Gamma $ and $\Phi\circ d\circ \mathrm{R}_{D^f}$.
    This shows the commutativity of the third square and completes the proof.
\end{proof}

We now note some corollaries of this theorem.

\begin{cor}
    Let $M$ be a closed oriented manifold and let $N$ be a closed manifold.
    Assume that $i_1,i_2\colon N\to M$ are two homotopic embeddings of $N$ into $M$.
    Then the homologies $\mathrm{H}_{\bullet}(P_N^{i_1} M)$ and $\mathrm{H}_{\bullet}(P_N^{i_2}M)$ are isomorphic as modules over the Chas-Sullivan ring of $M$. 
\end{cor}

\begin{cor}
    Let $M$ be a closed oriented manifold and let $i\colon N\to M$ be a closed submanifold with null-homotopic inclusion.
    Then the module structure on $\mathrm{H}_{\bullet}(P_N M)$ over the Chas-Sullivan ring of $M$ is isomorphic to the module structure $\nu_{N,\Omega}$ as described in Section \ref{sec_module_structure_explicit}.
\end{cor}

This last corollary gives us a very explicit way of computing the module structure in cases where the inclusion of the submanifold is null-homotopic.
Using the invariance of both the path product as well as the module structure combined with the results of Propositions \ref{prop_trivial_map_path_space} and \ref{prop_trivial_map_module} we shall now study the interaction between the path product and the module structure in this situation.

Recall that the Pontryagin ring $(\mathrm{H}_{\bullet}(\Omega M),\star)$ is in general not graded commutative.
For example if $M  = \mathbb{S}^{2n}$ then we have
$$  \mathrm{H}_{\bullet}(\Omega \mathbb{S}^n,\mathbb{Q}) \cong \mathbb{Q}[a] \qquad \text{with}\quad |a| = 2n-1 .    $$
If the Pontryagin ring was graded commutative, we would have
$$  a^2  = - a^2 = 0,    $$
but this is obviously not satisfied.
By the \textit{center} of a graded ring $R$ we mean the elements $x\in R$ such that
$$  xy =  (-1)^{|x||y|} yx \quad \text{for all}\,\,\, y\in R .    $$
In the example of the sphere, the center is given by the subring $\mathbb{Q}[a^2]$.

\begin{theorem}\label{theorem_algebra_over_cs}
    Let $M$ and $N$ be closed oriented manifolds and let $f\colon N\to M$ be a null-homotopic map.
    Let $n = \mathrm{dim}(M)$ and $k = \mathrm{dim}(N)$.
    Take homology with coefficients in a field and assume that the image of the Gysin map $j_! \colon \mathrm{H}_i(\Lambda M) \to \mathrm{H}_{i-n}(\Omega M)$ lies in the center of the Pontryagin ring of $\Omega M$.    
    Then the homology $\mathrm{H}_{\bullet}(P^f)$ equipped with the path product is an algebra over the Chas-Sullivan ring in the sense that the identity
    $$      A \cdot (X\wedge Y) =  (A\cdot X) \wedge Y   =  (-1)^{(|A| + n)(|X| + k)}      X \wedge (A\cdot Y)   $$
    holds for all $A\in\mathrm{H}_{\bullet}(\Lambda M)$ and $X,Y\in\mathrm{H}_{\bullet}(P^f)$.
\end{theorem}
\begin{proof}
    By Theorems \ref{theorem_invariance_path_product} and \ref{theorem_invariance_module_structure} as well as Propositions \ref{prop_trivial_map_path_space} and \ref{prop_trivial_map_module} we can compute the product and the module structure using the product $\mu_{N,\Omega}$ and the module structure $\nu_{N,\Omega}$.

    Let $\Phi\colon P^f\to N\times N\times \Omega M$ be the homotopy equivalence as in Section \ref{sec_invariance_product}.
    Let $A\in\mathrm{H}_{\bullet}(\Lambda M)$ and $X,Y\in\mathrm{H}_{\bullet}(P^f)$ such that
    $$   \Phi_*X =  a\otimes b\otimes x \quad \text{and}\quad \Phi_* Y =  c\otimes d\otimes y      $$
    for $a,b,c,d\in\mathrm{H}_{\bullet}(N)$ and $x,y\in\mathrm{H}_{\bullet}(\Omega M)$.
    Using the explicit formulas for $\mu_{N,\Omega}$ and $\nu_{N,\Omega}$ we compute
    \begin{eqnarray}\label{eq_comp_1}
            &  &  \nonumber \Phi_*(A\cdot (X\wedge Y)) = \\ &  &  (-1)^{|x|(|c|+|d|+k)}(-1)^{(|A|+n)(|a|+|d|+n)} \mu_{\beta}(a\otimes b\otimes c\otimes d\otimes d) \otimes (j_! A\star x\star y) .      
    \end{eqnarray}
    Furthermore, we get
    \begin{eqnarray}
        \label{eq_comp_2}
        &  &  \nonumber \Phi_*(A\cdot X)\wedge Y) =  \\  &  & (-1)^{(|A|+n)(|a|+|b|+|c|+|d|+n+k)}(-1)^{|x|(|c|+|d|+k)} \mu_{\beta}(a\otimes b\otimes c\otimes d) \otimes (j_! A\star x\star y) .
    \end{eqnarray}
    Using the fact that the expression $\mu_{\beta}(a\otimes b\otimes c\otimes d)$ is only non-zero if $|b|+|d| = k$ we see that the right hand sides of equations \eqref{eq_comp_1} and \eqref{eq_comp_2} agree.
    Finally, we compute
    \begin{eqnarray}\label{eq_comp_3}
         &  &  \nonumber \Phi_*(X \wedge (A\cdot Y)) =    \\  &  &  \nonumber (-1)^{(|A|+n)(|c|+|d|+n)}(-1)^{|x|(|c|+|d|+k)} \mu_{\beta}(a\otimes b\otimes c\otimes d) \otimes ( x \star j_! A\star y)   = \\ &  &
         (-1)^{(|A|+n)(|c|+|d|+n)}(-1)^{|x|(|c|+|d|+k)}(-1)^{(|A|+n)|x|} \mu_{\beta}(a\otimes b\otimes c\otimes d) \otimes (j_! A\star x\star y) 
    \end{eqnarray}
    where we used the assumption that $j_! A$ lies in the center of the Pontryagin ring in the second equality.
    A straight-forward computation yields that the right hand sides in equations \eqref{eq_comp_1} and \eqref{eq_comp_3} differ by the sign $(-1)^{(|A| + n)(|X| + k)}$.
    This completes the proof.
\end{proof}

\begin{remark}
    Let $M$ and $N$ be closed oriented manifolds and $f\colon N\to M$ a null-homotopic map.
    \begin{enumerate}
        \item Note that one might have expected the sign in the statement of Theorem \ref{theorem_algebra_over_cs} to be $(-1)^{|A||X|}$.
        However, if we regrade the homology of $\Lambda M$ by $-n$ and the homology of $P^f$ by $-k$ then we indeed get this "expected" sign.
        When considering the Chas-Sullivan product it is quite common to regrade the homology of the free loop space in exactly this way.
        \item 
    We note that the condition in the above theorem on the image of the Gysin map $j_!$ is clearly satisfied if the Pontryagin ring itself is graded commutative.
    It is well-known that this is the case for compact Lie groups and odd-dimensional spheres.
    Below we shall see that the condition is also satisfied for even-dimensional spheres.
    In fact, we are not aware of any example where this condition on the Gysin map is not satisfied.
    \end{enumerate}
\end{remark}

\begin{remark}
    Let $M$ and $N$ be closed manifolds and $f\colon N\to M$ a smooth map, not necessarily null-homotopic.
    It is a natural question whether the conclusion of Theorem \ref{theorem_algebra_over_cs} holds in this more general situation, i.e. is the homology $\mathrm{H}_{\bullet}(P^f)$ equipped with the path product an algebra over the Chas-Sullivan ring of $M$?
    The author conjectures that the answer to this question is negative in general.
    The reason for this is the following heuristic observation.
    Let $A\in \mathrm{H}_{\bullet}(\Lambda M) $ and $X,Y\in\mathrm{H}_{\bullet}(P^f)$.
    If we take representing cycles of these homology classes then recall that for the pairing of $A$ and $X$ the module structure intersects the basepoints of the loops in $A$ with the starting points of the paths in $X$.
    The path product of $X$ and $Y$ intersects the endpoints of the paths in $X$ with the starting points of the paths in $Y$.
    Assume that $A\cdot (X\wedge Y)$ is non-trivial, then these respective intersections of basepoints are non-trivial.
    If one finds a situation where the basepoints of $A$ and the starting points of $Y$ intersect trivially, then one would have $A\cdot (X\wedge Y) \neq 0$ while $X\wedge (A\cdot Y) = 0$.
    In this case the homology $\mathrm{H}_{\bullet}(P^f)$ equipped with the path product would not be an algebra over the Chas-Sullivan ring of $M$.
    It would be desirable to find a situation where this intuition could be made precise.
\end{remark}

We conclude this section by considering some specific situations.
First, we show that the conclusion of Theorem \ref{theorem_algebra_over_cs} always hold for spheres.
For odd-dimensional spheres we have already seen in the above remark that the condition on the Gysin map $\mathrm{H}_i(\Lambda \mathbb{S}^{2n+1})\to \mathrm{H}_{i-{2n}-1}(\Omega \mathbb{S}^{2n+1})$ to lie in the center of the Pontryagin ring is satisfied, since the Pontryagin ring of $\Omega \mathbb{S}^{2n+1}$ is graded commutative.
We now show that the image of the Gysin map $j_!$ is equal to the center of the Pontryagin ring for even-dimensional spheres.
\begin{lemma}\label{lemma_even_dim_sphere_gysin}
    Let $M = \mathbb{S}^{2n}$ be an even-dimensional sphere and take homology with coefficients in a field of characteristic not equal to $2$.
    Then the image of the Gysin map
    $$j_!\colon \mathrm{H}_i(\Lambda \mathbb{S}^{2n}) \to \mathrm{H}_{i-{2n}}(\Omega\mathbb{S}^{2n}) \cong \mathbb{K}[a]$$ 
    with $|a| = 2n-1$ is equal to the center $\mathbb{K}[a^2]$ of the Pontryagin ring.
\end{lemma}
\begin{proof}
     Consider the free loop fibration $\Lambda \mathbb{S}^{2n}\to \mathbb{S}^{2n}$ and let $\{E^m,\mathrm{d}_m\}_{m\in\mathbb{N}}$ be the associated Serre spectral sequence.
    Let $\{\widetilde{E}^m,\mathrm{d}_m\}_{m\in\mathbb{N}}$ be the Serre spectral sequence associated to the trivial fibration $\Omega\mathbb{S}^n\to \{\mathrm{pt}\}$.
    Of course the spectral sequence $\widetilde{E}^m$ is trivial.
     By \cite[Proposition 4.3]{meier:2011} the Gysin map $$j_!\colon \mathrm{H}_{\bullet}(\Lambda\mathbb{S}^{2n})\to \mathrm{H}_{\bullet-n}(\Omega \mathbb{S}^{2n})$$ induces a morphism of spectral sequences $F^m(j_!)\colon E^m\to \widetilde{E}^m$ of bidegree $(-n,0)$ such that on the $E_2$-page we have that
    $$    F^2(j_!)\colon E^2_{p,q}  \cong \mathrm{H}_p(\mathbb{S}^{2n})\otimes \mathrm{H}_q(\Omega \mathbb{S}^{2n}) \, \to \,\widetilde{E}^2_{p-n,q} \cong \mathrm{H}_{p-n}(\{\mathrm{pt} \}) \otimes \mathrm{H}_q(\Omega \mathbb{S}^{2n})   $$
    is given by $ F^2(j_!) = i_!\otimes \mathrm{id}$  
    where $i\colon \{\mathrm{pt}\}\hookrightarrow \mathbb{S}^{2n}$ is the inclusion of the basepoint.
    By the remark after \cite[Proposition 4.3]{meier:2011} this map of spectral sequences converges to the map induced by the Gysin morphism $j_!$ if $p$ is a smooth fiber bundle of Hilbert manifolds.
    This is indeed the case by \cite[Theorem 1.1]{chataur:2015}.
        The differential $\mathrm{d}_n$ on the $E^n$-page behaves as follows, see \cite{cohen:2003}.
    We have that $$   \mathrm{d}_n\colon E^n_{n,(2l-1)(n-1)}\to E^n_{0,2l(n-1)} \quad \text{for} \,\,\, l\in\mathbb{N}$$ is multiplication with $2$ under appropriate identifications of $E^n_{n,l(n-1)}$ and $ E^n_{0,(l+1)(n-1)}$ with $\mathbb{K}$.
    Therefore this differential is an isomorphism while all other differentials vanish.
    Therefore we get $$E^{\infty}_{n,(2l-1)(n-1)}\cong  E^{n+1}_{n,(2l-1)(n-1)} \cong \{0\}  $$
    for all $l\in\mathbb{N}$.
    Consequently the map $$F^{n+1}(j_!) \colon E^{n+1}_{n,(2l-1)(n-1)}  \to \widetilde{E}^{n+1}_{0,(2l-1)(n-1)}$$ vanishes for all $l\in\mathbb{N}$.
    Note that we have
    $$  \widetilde{E}^{n+1}_{0,(2l-1)(n-1)} \cong \mathrm{H}_{(2l-1)(n-1)}(\Omega \mathbb{S}^{2n})\cong  \mathbb{K}\langle a^{2l-1} \rangle . $$
    Thus we see that the image of $F^{n+1}(j_!)$ is equal to $\bigoplus_{i=0}^{\infty} \mathbb{K}\langle a^{2l}\rangle$ which is precisely the center of the Pontryagin ring of $\mathbb{S}^{2n}$.
    This shows the claim. 
\end{proof}

Since the inclusion of any submanifold into a sphere is null-homotopic we obtain the following corollary.
\begin{cor}
    Let $M = \mathbb{S}^n$ be a sphere and let $N\hookrightarrow \mathbb{S}^n$ be a closed oriented submanifold.
    Take homology with coefficients in a field.
    Then the homology $\mathrm{H}_{\bullet}(P_N\mathbb{S}^n)$ equipped with the path product is an algebra over the Chas-Sullivan ring of $\mathbb{S}^n$.
\end{cor}

We end this section by showing that in the case of Lie groups and spheres the homology $\mathrm{H}_{\bullet}(P^f)$ is a finitely generated module over the Chas-Sullivan ring.

\begin{prop}
    Let $M$ be a compact Lie group or a sphere and let $N$ be a closed manifold with a smooth map $f\colon N\to M$.
    Take homology with coefficients in a field $\mathbb{K}$ and assume that $f$ is null-homotopic.
    Then the homology $\mathrm{H}_{\bullet}(P^f)$ is finitely generated as a module over the Chas-Sullivan ring of $M$.
\end{prop}
\begin{proof}
     We begin with the case of $M= G$ being a compact Lie group.
     Note that the free loop space of a Lie group is homeomorphic to the product 
     $$     \Lambda G\cong G\times \Omega_e G   $$
     under the map $\gamma \mapsto (\gamma(0),\gamma(0)^{-1}\gamma)$.
     Thus the inclusion $j\colon \Omega_e G\to \Lambda G$ is just the inclusion of the second factor.
     Moreover, we have an isomorphism $\mathrm{H}_{\bullet}(\Lambda G)\cong \mathrm{H}_{\bullet}(G)\otimes \mathrm{H}_{\bullet}(\Omega G)$.
     It is easy to see that the Gysin map $j_!$ behaves as follows under this isomoprhism.
     We have
     $$   j_! ([G]\otimes X) = X     $$
     where $[G]\in\mathrm{H}_{n}(G)$ is the fundamental class and $X\in\mathrm{H}_{\bullet}(\Omega G)$ is arbitrary.
     If $A \in\mathrm{H}_{i}(G)$ with $i< n$ then we have
     $$    j_! (A\otimes X) = 0 .   $$
     Now, assume that $N$ is a closed manifold with a smooth null-homotopic map $f\colon N\to G$.
    Let $x_1,\ldots,x_m\in\mathrm{H}_{\bullet}(N)$ be a family of homology classes which generate $\mathrm{H}_{\bullet}(N)$ additively.
    Recall that we have an isomorphism $\mathrm{H}_{\bullet}(P^f) \cong \mathrm{H}_{\bullet}(N)\otimes \mathrm{H}_{\bullet}(N)\otimes\mathrm{H}_{\bullet}(\Omega G)$.
    We claim that the finite family of homology classes $S =\{x_i\otimes x_j \otimes \mathbb{1}\}_{i,j\in\{1,\ldots,m\}}$ is a generating set for $\mathrm{H}_{\bullet}(P^f)$ where $\mathbb{1}\in\mathrm{H}_0(\Omega_e G)$ is the unit of the Pontryagin ring.
    Let $X = x_i\otimes x_j\otimes b\in \mathrm{H}_{\bullet}(N)^{\otimes 2}\otimes \mathrm{H}_{\bullet}(\Omega G)$.
    Define $X_0 = x_i\otimes x_j \otimes \mathbb{1}$.
    Then we have
    \begin{equation}\label{eq_xi_lie_groups}
           \nu_{N,\Omega}( [G]\otimes b \otimes x_i\otimes x_j\otimes \mathbb{1}) =\pm x_i\otimes x_j\otimes b .  
    \end{equation}
    This shows that the set $S$ is a finite generating set of $\mathrm{H}_{\bullet}(P_N G)$ as a module over the Chas-Sullivan ring.

    For the second part note that with similar arguments as in the proof of Lemma \ref{lemma_even_dim_sphere_gysin} one can show that the image of the Gysin map $j_!\colon \mathrm{H}_{\bullet}(\Lambda \mathbb{S}^{2n+1})\to \mathrm{H}_{\bullet-n}(\Omega \mathbb{S}^{2n+1})$ is all of $\mathrm{H}_{\bullet}(\Omega \mathbb{S}^{2n+1})$.
    Thus for the proof of the proposition in the case of odd-dimensional spheres and for even-dimensional spheres with $\mathbb{Z}_2$-coefficients we can argue in the same way as we did for Lie groups.
    
    Finally, let $M = \mathbb{S}^{2n}$ be a sphere of even dimension and assume that $\mathbb{K}$ is a field with $\mathrm{char}(\mathbb{K})\neq 2$.
    Let $N$ be a closed manifold and let $f\colon N\to \mathbb{S}^{2n}$ be a smooth null-homotopic map.
    Let $x_1,\ldots, x_m\in\mathrm{H}_{\bullet}(N)$ be a basis for $\mathrm{H}_{\bullet}(N)$.
    Then for $l=0,1$ define the family
    $$    S_l  =   \{ x_i\otimes x_j \otimes a^l\}_{i,j\in\{1,\ldots,m\}}   $$
    where $a\in \mathrm{H}_{2n-1}(\Omega \mathbb{S}^{2n})$ is the generator of the Pontryagin ring.
    Set $S = S_0\cup S_1$.
    Again with similar arguments as above and using Lemma \ref{lemma_even_dim_sphere_gysin} one finds that $S$ is a generating set for $\mathrm{H}_{\bullet}(P^f)$ as a module over the Chas-Sullivan ring of $\mathbb{S}^{2n}$.   
\end{proof}

\bibliography{lit}
 \bibliographystyle{amsalpha}

\end{document}